\documentclass[11pt]{article}
\usepackage[left=1in,right=1in,top=1in,bottom=1in]{geometry}
\usepackage{times}
\usepackage{expl3}
\usepackage{cite}
\usepackage[table]{xcolor}
\usepackage{multirow}
\usepackage{stackengine} 
\usepackage{hhline}
\usepackage{lipsum}
\usepackage{titlesec}
\usepackage{wrapfig}
\usepackage{epsfig}
\usepackage{enumitem}
\usepackage{bbm}
\usepackage{calc}
\usepackage{graphicx}
\usepackage{amsmath}
\usepackage[title]{appendix}
\usepackage{amssymb}
\usepackage{epstopdf}
\usepackage{boldline}
\usepackage{calligra}
\usepackage{bm}
\usepackage{url}
\usepackage{blindtext}
\usepackage{accents}
\newcommand{\dbtilde}[1]{\accentset{\approx}{#1}}

\newcommand{\define}{\stackrel{\mbox{\tiny def}}{=}}

\newtheorem{definition}{Definition}
\newtheorem{theorem}{Theorem}
\newtheorem{corollary}{Corollary}
\newtheorem{lemma}{Lemma}

\usepackage{mathtools}
\usepackage{epstopdf}
\usepackage{balance}
\usepackage{thmtools}
\usepackage{thm-restate}
\usepackage{hyperref}
\usepackage{cleveref}

\usepackage[ruled,vlined]{algorithm2e}
\include{pythonlisting}
\newcommand{\ostar}{\mathbin{\mathpalette\make@circled\star}}

\makeatletter
\newcommand{\removelatexerror}{\let\@latex@error\@gobble}
\makeatother
\makeatletter
\newcommand*{\rom}[1]{\expandafter\@slowromancap\romannumeral #1@}
\makeatother

\ExplSyntaxOn
\newcommand\latinabbrev[1]{
  \peek_meaning:NTF . {
    #1\@}%
  { \peek_catcode:NTF a {
      #1.\@ }%
    {#1.\@}}}
\ExplSyntaxOff

\titleclass{\subsubsubsection}{straight}[\subsubsection]

\begin{document}
\vspace{1cm}
\title{Random Tensor Inequalities and Tail bounds for Bivariate Random Tensor Means, Part II}\vspace{1.8cm}
\author{Shih~Yu~Chang 
\thanks{Shih Yu Chang is with the Department of Applied Data Science,
San Jose State University, San Jose, CA, U. S. A. (e-mail: {\tt
shihyu.chang@sjsu.edu}).
           }}

\maketitle

\begin{abstract}
This is Part II of our work about random tensor inequalities and tail bounds for bivariate random tensor means. After reviewing basic facts about random tensors, we first consider tail bounds with more general connection functions. Then, a general Lie-Trotter formula for tensors is derived and this formula is applied to establish tail bounds for bivariate random tensor means involving tensor logarithm. All random tensors studied in our Part I work are assumed as positive definite (PD) random tensors, which are invertible tensors. In this Part II work, we generalize our tail bounds for bivariate random tensor means from positive definite (PD) random tensors to positive semidefinite (PSD) random tensors by defining Random Tensor Topology (RTT) and developing the limitation method based on RTT. Finally,  we apply our theory to establish tail bounds and  L\"owner ordering relationships for bivariate random tensor means before and after two tensor data processing methods: data fusion and linear transform. 
\end{abstract}

\begin{keywords}
Lie-Trotter formula, random tensors, bivariate tensor mean, L\"owner ordering, 
majorization ordering, tensor data processing.
\end{keywords}

\section{Introduction}

Ando-Hiai type inequalities for operator geometric means was first proved at a noble work in~\cite{ando1994log}. Since then, these inequalities have initiated active research in operator theory, e.g., general operator means, multivariable operator means, quantum information, etc~\cite{hiai2019log,hiai2017different}. When $g$ is a connection function associated to a particular operator mean and two positive invertible operators $\bm{A}$ and $\bm{B}$, the Ando-Hiai type inequalities can be expressed as:
\begin{eqnarray}
\bm{A} \#_g \bm{B} &\preceq& \bm{I} \implies \bm{A}^q \#_g \bm{B}^q \preceq \bm{I},\\
\bm{A} \#_g \bm{B} &\succeq& \bm{I} \implies \bm{A}^q \#_g \bm{B}^q \succeq \bm{I},
\end{eqnarray}
where $q \geq 1$. In the Part I of our work, we discuss tail bounds for bivariate random tensor means for various types of connection functions $g$ and broader range of exponent $q$. In this paper, we consider tail bounds with more general connection function $G_n(x)$, where $G_n(x)$ is defined by $G_n(x) \define x^n g(x)$ for $n \in \mathbb{N}$ by considering operators as random tensors. There are many works dedicated to topics about random tensors recently~\cite{chang2022TWF,chang2022randouble,chang2022tPII,chang2022tPI,
chang2022generaltail,chang2022TKF,chang2022tensorq}.  

The Lie-Trotter formula is a fundamental result in numerical analysis that allows the approximation of solutions to differential equations using a sequence of simpler sub-problems. Specifically, the formula states that the solution to a composite system of differential equations can be approximated by applying the solutions to the individual sub-problems in a sequential manner. The Lie-Trotter formula has applications in a wide range of fields. Several examples of Lie-Trotter formula applications are discussed as follows. Computational physics: The Lie-Trotter formula is used to simulate the behavior of physical systems that can be modeled by differential equations, such as the Schrödinger equation in quantum mechanics or the Navier-Stokes equation in fluid dynamics. By approximating the solutions to these equations using the Lie-Trotter formula, investigators can gain insights into the behavior of complex systems and make predictions about their behavior under different conditions. Molecular dynamics: The Lie-Trotter formula is used to simulate the behavior of molecular systems, such as proteins or nucleic acids, by approximating the solutions to the equations of motion for each atom or molecule in the system. This allows investigators to study the behavior of these complex systems and to design new drugs or materials based on their properties. Control theory: The Lie-Trotter formula is used in the design of control systems for complex engineering applications, such as aerospace or robotics. By approximating the solutions to the equations governing the behavior of these systems, engineers can design control systems that can stabilize and optimize the behavior of these systems in real-time. Financial mathematics: The Lie-Trotter formula is used in the pricing of financial derivatives, such as options or futures contracts, by approximating the solutions to the partial differential equations governing the behavior of these financial instruments. This allows traders and investors to make informed decisions about their investments based on the predicted behavior of these instruments under different market conditions. In general, the Lie-Trotter formula is a powerful tool that enables investigators and engineers to simulate the behavior of complex systems and to make informed decisions based on the predicted behavior of these systems under different conditions~\cite{ahn2007extended}. In this paper, we will apply the Lie-Trotter formula to establish tail bounds for bivariate random tensor means involving tensor logarithm.

All random tensors discussed in our Part I work are assumed as positive definite (PD) random tensors, which are invertible tensors. In this work, we will generalize our tail bounds for bivariate random tensor means from positive definite (PD) random tensors to positive semidefinite (PSD) random tensors. We first define the notion about Random Tensor Topology (RTT) and develop a limitation method based on RTT to generalize tail bounds for bivariate random tensor means from random positive definite (PD) tensors to random positive semidefinite (PSD) tensors. The derivation of the Lie-Trotter formula and the aforementioned limitation method is derived based on the work~\cite{hiai2020ando}.

Graph data fusion is the process of combining information from multiple graphs or network datasets into a single, unified graph. The goal of graph data fusion is to leverage the strengths of each individual graph to create a more comprehensive and accurate representation of the underlying data. There are several techniques used in graph data fusion. First is graph merging. This technique involves combining the nodes and edges of two or more graphs into a single graph. This can be done by identifying common nodes or edges and merging them, or by creating new nodes and edges to connect the different graphs. The second is graph alignment. Graph alignment involves finding correspondences between nodes or edges in different graphs and aligning them. This can be done by comparing node or edge attributes, or by using graph matching algorithms. Finally, graph coarsening involves reducing the complexity of the graphs by grouping nodes together into clusters or communities. This can be done by using clustering algorithms or community detection algorithms. Graph data fusion has applications in a variety of fields, including social network analysis, bioinformatics, and transportation planning. It allows investigators and analysts to gain a more complete understanding of complex systems and make more accurate predictions based on the data~\cite{cichocki2014era}. As tensors (matrices) are ideal tools to perform algebraic operations over graph data, in this work, we apply our theory about tensor convexity to develop tensor mean fusion inequality. 

Another important data processing method applied to graph data is linear transform over graph data. Graph data processing by linear transform refers to a mathematical operation applied to graphs that transforms the original graph into a new graph with different properties while preserving certain features of the original graph. In particular, a linear transform of a graph is a function that maps the graph's adjacency matrix (tensor) to a new matrix (tensor) using a linear operator. This new matrix represents a transformed version of the original graph. Some common examples of linear transforms are discussed below. Laplacian transform: This transform is based on the graph Laplacian matrix and is commonly used in spectral graph theory. It can be used to identify clusters or communities in a graph and to detect graph properties such as connectivity and symmetry. Graph Fourier transform: This transform is based on the graph's eigenvalues and eigenvectors and is used to analyze the frequency content of a graph. It can be used to identify patterns or motifs in a graph and to compare different graphs based on their spectral properties. Wavelet transform: This transform is used to analyze the local properties of a graph by decomposing the graph into different scales and analyzing each scale separately. It can be used to identify local patterns or anomalies in a graph and to filter out noise or irrelevant information. Linear transforms of graphs have many applications in data processing, including image processing, signal processing, and machine learning. They can be used to extract meaningful features from graphs, to reduce the dimensionality of graph data, and to enable efficient computation on large graphs~\cite{zhang2015graph}. We apply operator (tensor) Jensen inequality to prove tensor mean inequality for bivariate tensor means under linear transform. 

The organization of this Part II work is summarized as follows. Those important concepts like connection functions classifications, Kantorovich contant, and random tensor inequalities discussed in the Part I of this work are restated here in Section~\ref{sec:Preliminaries} for self-contained presentation. In Section~\ref{sec:Tail Bounds for Bivariate Random Tensor Means of Generalized Connection Functions}, we consider tail bounds with more general connection function $F_n(x)$, where $F_n(x)$ is defined by $F_n(x) \define x^n f(x)$. In Section~\ref{sec:Lie-Trotter Formula and Its Applications}, a general Lie-Trotter formula for tensors is derived and this formula is applied to establish tail bounds for bivariate random tensor means involving tensor logarithm. In section~\ref{sec:Loewner and Majorization Ordering for Bivariate Random Tensor Means with Non-invertible Tensors}, a new topology to characterize the limitation behavior of random tensor, namely Random Tensor Topology (RTT), is first introduced, then we develop a limitation method based on RTT to generalize tail bounds for bivariate random tensor means from random positive definite (PD) tensors to random positive semidefinite (PSD) tensors. Finally, we apply our theory about tail bounds for bivariate random tensors to tensor data processing by characterizing L\"owner orders and tail bounds relationships for bivariate random tensor means after transformations (data processing).

\section{Preliminaries}\label{sec:Preliminaries}

In this section, we will review required definitions and facts that will be used in later sections. Let us represent the Hermitian eigenvalues of a Hermitian tensor $\mathcal{H} \in \mathbb{C}^{I_1 \times \dots \times I_N \times I_1 \times \dots \times I_N} $ in decreasing order by the vector $\vec{\lambda}(\mathcal{H}) = (\lambda_1(\mathcal{H}), \cdots, \lambda_r(\mathcal{H}))$, where $r$ is the Hermitian rank of the tensor $\mathcal{H}$. We use $\mathbb{R}_{\geq 0} (\mathbb{R}_{> 0})$ to represent a set of nonnegative (positive) real numbers. Let $\left\Vert \cdot \right\Vert_{\rho}$ be a unitarily invariant tensor norm, i.e., $\left\Vert \mathcal{H}\star_N \mathcal{U}\right\Vert_{\rho} = \left\Vert \mathcal{U}\star_N \mathcal{H}\right\Vert_{\rho} = \left\Vert \mathcal{H}\right\Vert_{\rho} $,  where $\mathcal{U}$ is for unitary tensor. Let $\rho : \mathbb{R}_{\geq 0}^r \rightarrow \mathbb{R}_{\geq 0}$ be the corresponding gauge function that satisfies H${\rm \ddot{o}}$lder’s inequality so that 
\begin{eqnarray}\label{eq:def gauge func and general unitarily invariant norm}
\left\Vert \mathcal{H} \right\Vert_{\rho} = \left\Vert |\mathcal{H}| \right\Vert_{\rho} = \rho(\vec{\lambda}( | \mathcal{H} | ) ),
\end{eqnarray}
where $ |\mathcal{H}|  \define \sqrt{\mathcal{H}^H \star_N \mathcal{H}} $.

We define following sets of positive functions.
\begin{eqnarray}\label{eq:three sets of functions}
\mbox{TMI}&=&\{f: \mbox{tensor monotone increasing on $(0, \infty)$, $f $$>$$ 0$}\};\nonumber \\
\mbox{TMD}&=&\{g: \mbox{tensor monotone decreasing on $(0, \infty)$, $g $$>$$ 0$}\};\nonumber \\
\mbox{TC}&=&\{h: \mbox{tensor convex on $(0, \infty)$, $h $$>$$ 0$}\}.
\end{eqnarray}
\begin{eqnarray}\label{eq:three sets of functions with one}
\mbox{TMI}^{1}&=&\{f: \mbox{tensor monotone increasing on $(0, \infty)$, $f $$>$$ 0$ and $f(1) $$=$$ 1$}\};\nonumber \\
\mbox{TMD}^{1}&=&\{g: \mbox{tensor monotone decreasing on $(0, \infty)$, $g $$>$$ 0$ and $g(1) $$=$$ 1$}\};\nonumber \\
\mbox{TC}^{1}&=&\{h: \mbox{tensor convex on $(0, \infty)$, $h $$>$$ 0$ and $h(1) $$=$$ 1$}\}.
\end{eqnarray}

The default product between two tensors is $\star_N$. We will specfiy the exact product symbol if the product is not $\star_N$.

Lemmas~\ref{lma:2.3} and~\ref{lma:Kantorovich type inequality} about Kantorovich type inequality is proved in our Part I work. 
\begin{lemma}\label{lma:2.3}
Let $\mathcal{A}, \mathcal{B} \in \mathbb{C}^{I_1 \times \cdots \times I_N \times I_1 \times \cdots \times I_N}$ be two SPD tensors with $\mathcal{A} \succeq \mathcal{B}$,  and let $q $$\in$$ [0,1]$. Then, we have
\begin{eqnarray}\label{eq1:lma:2.3}
\mathcal{A}^q \succeq \mathcal{B}^q. 
\end{eqnarray}
\end{lemma}

\begin{lemma}\label{lma:Kantorovich type inequality}
Let $\mathcal{A}, \mathcal{B} $$\in$$ \mathbb{C}^{I_1 \times \cdots \times I_N \times I_1 \times \cdots \times I_N}$ be two PD tensors such that 
\begin{eqnarray}\label{eq1:lma:Kantorovich type inequality}
m_1 \mathcal{I} \preceq \mathcal{A} \preceq M_1 \mathcal{I},~\mbox{and} \nonumber \\
m_2 \mathcal{I} \preceq \mathcal{B} \preceq M_2 \mathcal{I},~~~~~~~
\end{eqnarray}
where $M_1 $$>$$ m_1 $$>$$ 0$, and $M_2 $$>$$ m_2 $$>$$ 0$. If $\mathcal{B} $$\preceq$$ \mathcal{A}$, we have
\begin{eqnarray}\label{eq2:lma:Kantorovich type inequality}
\mathcal{B}^p &\preceq& \mathrm{K}(m_1,M_1,p) \mathcal{A}^p, \nonumber \\
\mathcal{B}^p &\preceq& \mathrm{K}(m_2,M_2,p) \mathcal{A}^p,
\end{eqnarray}
where the Kantorovich contant, $\mathrm{K}(m,M,p)$, can be expressed by
\begin{eqnarray}
\mathrm{K}(m,M,p) &=& \left(\frac{(p-1)\left(M^p - m^p\right)}{p\left(mM^p - Mm^p\right)}\right)^p\frac{mM^p - Mm^p}{(p-1)(M-m)}. 
\end{eqnarray}
\end{lemma}

Lemma~\ref{lma:Loewner ordering with Markov Cheb inequalities} below provides tail bounds according to 
L\"owner ordering between random tensors. 
\begin{lemma}\label{lma:Loewner ordering with Markov Cheb inequalities}
Given the following random PD tensors $\mathcal{X}, \mathcal{Y}, \mathcal{Z} $$\in$$ \mathbb{C}^{I_1 \times \cdots \times I_N \times I_1 \times \cdots \times I_N}$ with the relation $\mathcal{X} $$\preceq$$ \mathcal{Y} $$\preceq$$ \mathcal{Z}$, and a deterministic PD tensor $\mathcal{C}$, we have 
\begin{eqnarray}\label{eq1:lma:Loewner ordering with Markov Cheb inequalities}
\mathrm{Pr}\left(\mathcal{Y} \npreceq \mathcal{C}\right) &\leq& \mathrm{Tr}\left(\mathbb{E}\left[\mathcal{Z}^q\right] \star_N \mathcal{C}^{-1}\right),
\end{eqnarray}
where $q \geq 1$. We also have 
\begin{eqnarray}\label{eq2:lma:Loewner ordering with Markov Cheb inequalities}
\mathrm{Pr}\left(\mathcal{X} \npreceq \mathcal{C}\right) &\leq& \mathrm{Tr}\left(\mathbb{E}\left[\mathcal{Y}^q\right] \star_N \mathcal{C}^{-1}\right).
\end{eqnarray}
\end{lemma}

\section{Tail Bounds for Bivariate Random Tensor Means of Generalized Connection Functions}\label{sec:Tail Bounds for Bivariate Random Tensor Means of Generalized Connection Functions}

In this section, we will consider tail bounds for bivariate random tensor means of generalized connection functions. We define the following generalized connection function $F_n(x)$ as $F_n(x) \define x^n f(x)$, where $f(x)$ comes from $\mbox{TMI}^{1}$. The bivariate random tensor mean with respect to $F_n(x)$ can be expressed as 
\begin{eqnarray}\label{eq:F_n def}
\mathcal{X} \#_{F_n} \mathcal{Y} = \mathcal{X} \#_{x^n f(x)} \mathcal{Y}.
\end{eqnarray}
Then, given random tensors $\mathcal{X} \in \mathbb{C}^{I_1 \times \cdots \times I_N \times I_1 \times \cdots \times I_N}$ and $\mathcal{Y} \in \mathbb{C}^{I_1 \times \cdots \times I_N \times I_1 \times \cdots \times I_N}$, we have the following recursive relation:
\begin{eqnarray}\label{eq:recursive}
\mathcal{X} \#_{F_n} \mathcal{Y}=\mathcal{X} \star_N \left(\mathcal{Y}^{-1} \#_{F_{n-1}} \mathcal{X}^{-1}\right) \star_N \mathcal{X} = \mathcal{X} \star_N \mathcal{Y}^{-1} \star_N \left(\mathcal{X} \#_{F_{n-2}} \mathcal{Y}\right) \star_N \mathcal{Y}^{-1}\star_N \mathcal{X}
\end{eqnarray}
where $n \in \mathbb{N}$ and $n \geq 2$. We will present several Lemmas before establishing our main result in this section.  

\begin{lemma}\label{lma:inv func TMI lma 4.2}
Let $f \in \mbox{TMI}^{1}$, $n \in \mathbb{N}$ and $g(x) \define (x^n f(x))^{[-1]}$ be the inverse function of the function $x^n f(x)$ on $(0, \infty)$. Then, $g^m(x) \in \mbox{TMI}^{1}$ for any $m \in \{0,1,2,\cdots,n\}$. 
\end{lemma}
\textbf{Proof:}
For $n=1$, it is proved by Lemma 5 in~\cite{ando1988comparison} that $(xf(x))^{[-1]} \in \mbox{TMI}^{1}$ given that $f \in \mbox{TMI}^{1}$. If $n \geq 2$, we have
\begin{eqnarray}
g^{n}(x) &=& \left((x^n f(x))^{[-1]}\right)^n = \left(((x f_n(x))\circ x^n)^{[-1]}\right)^n \nonumber \\
&=& \left((x^n)^{[-1]}\circ (xf_n(x))^{[-1]}\right)^n  \nonumber \\
&=& (xf_n(x))^{[-1]},
\end{eqnarray}
where $f_n(x) \define f(x^{1/n})$. Because $f_n(x) \in \mbox{TMI}^{1}$, we have $(xf_n(x))^{[-1]}\in \mbox{TMI}^{1}$. Then, this Lemma is proved. 
$\hfill\Box$

Following two Lemmas are about the bivariate random tensor means for $\mathcal{X} \#_{F_{2m-1}} \mathcal{Y}$ (odd exponent) and  $\mathcal{X} \#_{F_{2m}} \mathcal{Y}$ (even exponent), where $m \in \mathbb{N}$. 
\begin{lemma}[odd exponent]\label{lma:odd 4.3}
Let $f \in \mbox{TMI}^{1}$ with pmi property and $m \in \mathbb{N}$ with $m \geq 2$. For  $q > 0$, we also assume that 
\begin{eqnarray}\label{eq0:lma:odd 4.3}
f(\mathcal{Z}^q)\preceq M_1 f^q(\mathcal{Z}),
\end{eqnarray}
where $\mathcal{Z}$ is an arbitrary PD tensor and $M_1 \geq 1$ is a positive constant. Besides, we define the function $F(x)$ as $F(x) \define x^{2m-2}f(x)$, the function $g(x)$ as $g(x) \define \left(F^{[-1]}(x^{-1}) \right)^{-1}$, and the function $h(x)$ as $h(x) \define x g(x)$. Given two PD random tensors $\mathcal{X} \in \mathbb{C}^{I_1 \times \cdots \times I_N \times I_1 \times \cdots \times I_N}$ and $\mathcal{Y} \in \mathbb{C}^{I_1 \times \cdots \times I_N \times I_1 \times \cdots \times I_N}$, if $\mathcal{X} \#_{x^{2m-1}f(x)} \mathcal{Y}\preceq \mathcal{I}$ almost surely and $q >0$, we have
\begin{eqnarray}\label{eq1:lma:odd 4.3}
\mathcal{X}^q \#_{x^{2m-1}f(x)} \mathcal{Y}^q \preceq M_1\left(\prod_{k=2}^{m}\mathrm{K}_k\right)\mathcal{I},
\end{eqnarray}
where $\mathrm{K}_k$ are Kantorovich constants defined as
\begin{eqnarray}\label{eq1-1:lma:odd 4.3}
\mathrm{K}_k \define \mathrm{K}(\lambda^{-1}_{\max}\left(\mathcal{X}^{-1}g^{(m-k)}(\mathcal{X})\right),\lambda^{-1}_{\min}\left(\mathcal{X}^{-1}g^{(m-k)}(\mathcal{X})\right),2q).
\end{eqnarray}
\end{lemma}
\textbf{Proof:}
Since all multiplications between tensors are $\star_N$, we will remove these notaions in this proof for simplification. From recursion relation given by Eq.~\eqref{eq:recursive}, we have 
\begin{eqnarray}
\mathcal{I} \succeq \mathcal{X}\#_{x^{2m-1}f(x)}\mathcal{Y} = \mathcal{X}\left(\mathcal{Y}^{-1}\#_{x^{2m-2}f(x)}\mathcal{X}^{-1}\right)\mathcal{X},
\end{eqnarray}
which is equivalent to 
\begin{eqnarray}\label{eq2:lma:odd 4.3}
\left(\mathcal{X}^{1/2}\mathcal{Y}^{-1}\mathcal{X}^{1/2}\right)^{2m-2}f\left(\mathcal{X}^{1/2}\mathcal{Y}^{-1}\mathcal{X}^{1/2}\right) \preceq \mathcal{X}^{-1}. 
\end{eqnarray}

From Eq.~\eqref{eq2:lma:odd 4.3}, we have 
\begin{eqnarray}\label{eq3:lma:odd 4.3}
\mathcal{X}^{-1/2}\mathcal{Y}\mathcal{X}^{-1/2} \succeq g(\mathcal{X})
\Longleftrightarrow
\mathcal{Y} \succeq h(\mathcal{X}).
\end{eqnarray}
From Lemma~\ref{lma:inv func TMI lma 4.2}, both functions $g(x)$ and $h^{[1]}(x)$ are in $\mbox{TMI}^{1}$.

For $n=1,2,\cdots,m$, we will prove the following inequality by induction.
\begin{eqnarray}\label{eq4:lma:odd 4.3}
\mathcal{X}^q \#_{x^{2n-1}f(x)} \mathcal{Y}^q \preceq M_1\left(\prod_{k=2}^{n}\mathrm{K}_k\right)g^{2(m-n)q}(\mathcal{X}),
\end{eqnarray}
where $\mathrm{K}_k$ is Kantorovich contants defined by Eq.~\eqref{eq1-1:lma:odd 4.3}.

For the case $n=1$, we have 
\begin{eqnarray}\label{eq4:lma:odd 4.3}
\mathcal{X}^q \#_{xf(x)} \mathcal{Y}^q &=& \mathcal{Y}^q \#_{f(1/x)} \mathcal{X}^q \nonumber \\
&\preceq_1& (h(\mathcal{X}))^q \#_{f(1/x)} \mathcal{X}^q \nonumber \\
&=& \mathcal{X}^q f\left(g^{-q}\left(\mathcal{X}\right)\right) \nonumber \\
&\preceq_2&  M_1\left(\mathcal{X}f\left(g^{-1}\left(\mathcal{X}\right)\right)\right)^q \nonumber \\
&=& M_1g^{2(m-1)q}(\mathcal{X})
\end{eqnarray}
where $\preceq_1$ comes from $\mathcal{Y} \succeq h(\mathcal{X})$, and $\preceq_2$ comes from pmi of $f$ for $0 < q \leq 1$ and comes from Eq.~\eqref{eq0:lma:odd 4.3} for $q \geq 1$. Suppose the inequality given by Eq.~\eqref{eq4:lma:odd 4.3} is true for $n \leq m-1$, we have 
\begin{eqnarray}
\mathcal{X}^q \#_{x^{2n+1}f(x)} \mathcal{Y}^q &=&
\mathcal{X}^q \mathcal{Y}^{-q}\left(\mathcal{X}^q \#_{x^{2n-1}f(x)} \mathcal{Y}^q \right)
\mathcal{Y}^{-q}\mathcal{X}^{q} \nonumber \\
&\preceq_1&\mathcal{X}^q \mathcal{Y}^{-q}\left( M_1\left(\prod_{k=2}^{n}\mathrm{K}_k\right)g^{2(m-n)q}(\mathcal{X})\right)
\mathcal{Y}^{-q}\mathcal{X}^{q} \nonumber \\
&\preceq_2&\mathcal{X}^q \mathcal{Y}^{-q}\left( M_1\left(\prod_{k=2}^{n}\mathrm{K}_k\right)g^{2(m-n)q}(h^{[-1]}(\mathcal{Y}))\right)
\mathcal{Y}^{-q}\mathcal{X}^{q}\nonumber \\
&=& M_1\left(\prod_{k=2}^{n}\mathrm{K}_k\right)\mathcal{X}^q\left((h^{[-1]}(\mathcal{Y}))^{-1}g^{(m-n-1)}(h^{[-1]}(\mathcal{Y}))\right)^{2q}\mathcal{X}^{q}\nonumber \\
&\preceq_3&  M_1\left(\prod_{k=2}^{n+1}\mathrm{K}_k\right)\mathcal{X}^q\left(\mathcal{X}^{-1}g^{(m-n-1)}(\mathcal{X})\right)^{2q}\mathcal{X}^{q} \nonumber \\
&=&M_1\left(\prod_{k=2}^{n+1}\mathrm{K}_k\right)g^{2(m-n-1)q}(\mathcal{X}),
\end{eqnarray}
where $\preceq_1$ comes from induction hypothesis, $\preceq_2$ comes from Lemma~\ref{lma:inv func TMI lma 4.2}, $\preceq_3$ comes from that the function $x/g^{m-n-1}(x) \in \mbox{TMI}^{1}$, the relation $\mathcal{Y} \succeq h(\mathcal{X})$ and Lemma~\ref{lma:Kantorovich type inequality} for $q \geq 1/2$. For  $0 < q \leq 1/2$, $\preceq_3$ comes from Lemma~\ref{lma:2.3}. This Lemma is proved by setting $n=m$.
$\hfill\Box$

\begin{lemma}[even exponent]\label{lma:even 4.4}
Let $f \in \mbox{TMI}^{1}$ with pmi property and $m \in \mathbb{N}$ with $m \geq 2$. For  $q > 0$, we also assume that 
\begin{eqnarray}\label{eq0:lma:even 4.4}
f(\mathcal{Z}^q)\preceq M_1 f^q(\mathcal{Z}),
\end{eqnarray}
where $\mathcal{Z}$ is an arbitrary PD tensor and $M_1 \geq 1$ is a positive constant. Besides, we define the function $F(x)$ as $F(x) \define x^{2m-1}f(x)$, the function $g(x)$ as $g(x) \define \left(F^{[-1]}(x^{-1}) \right)^{-1}$, and the function $h(x)$ as $h(x) \define x g(x)$. Given two PD random tensors $\mathcal{X} \in \mathbb{C}^{I_1 \times \cdots \times I_N \times I_1 \times \cdots \times I_N}$ and $\mathcal{Y} \in \mathbb{C}^{I_1 \times \cdots \times I_N \times I_1 \times \cdots \times I_N}$, if $\mathcal{X} \#_{x^{2m}f(x)} \mathcal{Y}\preceq \mathcal{I}$ almost surely and $q >0$, we have
\begin{eqnarray}\label{eq1:lma:even 4.4}
\mathcal{X}^q \#_{x^{2m}f(x)} \mathcal{Y}^q \preceq M_1\left(\prod_{k=1}^{m}\mathrm{K}_k\right)\mathcal{I},
\end{eqnarray}
where $\mathrm{K}_k$ are Kantorovich contants defined as
\begin{eqnarray}\label{eq1-1:lma:even 4.4}
\mathrm{K}_k \define \mathrm{K}(\lambda^{-1}_{\max}\left(\mathcal{X}^{-1}g^{(m-k)}(\mathcal{X})\right),\lambda^{-1}_{\min}\left(\mathcal{X}^{-1}g^{(m-k)}(\mathcal{X})\right),2q).
\end{eqnarray}
\end{lemma}
\textbf{Proof:}
From recursion relation given by Eq.~\eqref{eq:recursive}, we have 
\begin{eqnarray}
\mathcal{I} \succeq \mathcal{X}\#_{x^{2m}f(x)}\mathcal{Y} = \mathcal{X}\left(\mathcal{Y}^{-1}\#_{x^{2m-1}f(x)}\mathcal{X}^{-1}\right)\mathcal{X},
\end{eqnarray}
which is equivalent to 
\begin{eqnarray}\label{eq2:lma:even 4.4}
\left(\mathcal{X}^{1/2}\mathcal{Y}^{-1}\mathcal{X}^{1/2}\right)^{2m-1}f\left(\mathcal{X}^{1/2}\mathcal{Y}^{-1}\mathcal{X}^{1/2}\right) \preceq \mathcal{X}^{-1}. 
\end{eqnarray}

From Eq.~\eqref{eq2:lma:even 4.4}, we have 
\begin{eqnarray}\label{eq3:lma:even 4.4}
\mathcal{X}^{-1/2}\mathcal{Y}\mathcal{X}^{-1/2} \succeq g(\mathcal{X})
\Longleftrightarrow
\mathcal{Y} \succeq h(\mathcal{X}).
\end{eqnarray}
From Lemma~\ref{lma:inv func TMI lma 4.2}, both functions $g(x)$ and $h^{[1]}(x)$ are in $\mbox{TMI}^{1}$.

For $n=1,2,\cdots,m$, we will prove the following inequality by induction.
\begin{eqnarray}\label{eq4:lma:even 4.4}
\mathcal{X}^q \#_{x^{2n}f(x)} \mathcal{Y}^q \preceq M_1\left(\prod_{k=1}^{n}\mathrm{K}_k\right)g^{2(m-n)q}(\mathcal{X}),
\end{eqnarray}
where $\mathrm{K}_k$ is Kantorovich contants defined by Eq.~\eqref{eq1-1:lma:even 4.4}.

For the case $n=1$, we have 
\begin{eqnarray}\label{eq4:lma:even 4.4}
\mathcal{X}^q \#_{x^2f(x)} \mathcal{Y}^q &=& \mathcal{X}^q\mathcal{Y}^{-q}\left(\mathcal{Y}^{q}\#_{f(x)} \mathcal{X}^q\right)\mathcal{Y}^{-q}\mathcal{X}^q \nonumber \\
&\preceq_1& \mathcal{X}^q\mathcal{Y}^{-q}\left(\mathcal{Y}^{q}\#_{f(x)} \left(h^{[-1]}(\mathcal{Y})\right)^q\right)\mathcal{Y}^{-q}\mathcal{X}^q \nonumber \\
&=& \mathcal{X}^q\left(\mathcal{Y}^{-q}f\left(\left(\mathcal{Y}^{-1}h^{[-1]}(\mathcal{Y})\right)^q\right)\right)\mathcal{X}^q \nonumber \\
&\preceq_2&  M_1\mathcal{X}^{q}\left(\mathcal{Y}^{-1}f\left(\mathcal{Y}^{-1}h^{[-1]}(\mathcal{Y})\right)\right)^q\mathcal{X}^q \nonumber \\
&=& M_1\mathcal{X}^{q}\left((h^{[-1]}(\mathcal{Y}))^{-1}g^{m-1}\left(h^{[-1]}(\mathcal{Y})\right)\right)^{2q}\mathcal{X}^q \nonumber \\
&\preceq_3& M_1\mathrm{K}_1\mathcal{X}^{q}\left(\mathcal{X}^{-1}g^{m-1}\left(\mathcal{X}\right)\right)^{2q}\mathcal{X}^q \nonumber \\
&=& M_1\mathrm{K}_1g^{2(m-1)q}\left(\mathcal{X}\right),
\end{eqnarray}
where $\preceq_1$ comes from $\mathcal{Y} \succeq h(\mathcal{X})$, $\preceq_2$ comes from pmi of $f$ for $0 < q \leq 1$ and comes from Eq.~\eqref{eq0:lma:even 4.4} for $q \geq 1$, and $\preceq_3$ comes from that the function $x/g^{m-1}(x) \in \mbox{TMI}^{1}$, the relation $\mathcal{Y} \succeq h(\mathcal{X})$ and Lemma~\ref{lma:Kantorovich type inequality} for $q \geq 1/2$. For  $0 < q \leq 1/2$, $\preceq_3$ comes from Lemma~\ref{lma:2.3}.

Suppose the inequality given by Eq.~\eqref{eq4:lma:even 4.4} is true for $n \leq m-1$, we have 
\begin{eqnarray}
\mathcal{X}^q \#_{x^{2n+2}f(x)} \mathcal{Y}^q &=&
\mathcal{X}^q \mathcal{Y}^{-q}\left(\mathcal{X}^q \#_{x^{2n}f(x)} \mathcal{Y}^q \right)
\mathcal{Y}^{-q}\mathcal{X}^{q} \nonumber \\
&\preceq_1&\mathcal{X}^q \mathcal{Y}^{-q}\left( M_1\left(\prod_{k=1}^{n}\mathrm{K}_k\right)g^{2(m-n)q}(\mathcal{X})\right)
\mathcal{Y}^{-q}\mathcal{X}^{q} \nonumber \\
&\preceq_2&\mathcal{X}^q \mathcal{Y}^{-q}\left( M_1\left(\prod_{k=1}^{n}\mathrm{K}_k\right)g^{2(m-n)q}(h^{[-1]}(\mathcal{Y}))\right)
\mathcal{Y}^{-q}\mathcal{X}^{q}\nonumber \\
&=& M_1\left(\prod_{k=1}^{n}\mathrm{K}_k\right)\mathcal{X}^q\left((h^{[-1]}(\mathcal{Y}))^{-1}g^{(m-n-1)}(h^{[-1]}(\mathcal{Y}))\right)^{2q}\mathcal{X}^{q}\nonumber \\
&\preceq_3&  M_1\left(\prod_{k=1}^{n+1}\mathrm{K}_k\right)\mathcal{X}^q\left(\mathcal{X}^{-1}g^{(m-n-1)}(\mathcal{X})\right)^{2q}\mathcal{X}^{q} \nonumber \\
&=&M_1\left(\prod_{k=1}^{n+1}\mathrm{K}_k\right)g^{2(m-n-1)q}(\mathcal{X}),
\end{eqnarray}
where $\preceq_1$ comes from induction hypothesis, $\preceq_2$ comes from Lemma~\ref{lma:inv func TMI lma 4.2}, $\preceq_3$ comes from that the function $x/g^{m-n-1}(x) \in \mbox{TMI}^{1}$, the relation $\mathcal{Y} \succeq h(\mathcal{X})$ and Lemma~\ref{lma:Kantorovich type inequality} for $q \geq 1/2$. For  $0 < q \leq 1/2$, $\preceq_3$ comes from Lemma~\ref{lma:2.3}. This Lemma is proved by setting $n=m$.
$\hfill\Box$

We are ready to present the following main Theorem in this section.
\begin{theorem}\label{thm:4.1}
Let $f \in \mbox{TMI}^{1}$ with pmi property and $m \in \mathbb{N}$ with $m \geq 2$. For  $q > 0$, we also assume that 
\begin{eqnarray}\label{eq0:thm:4.1}
f(\mathcal{Z}^q)\preceq M_1 f^q(\mathcal{Z}),
\end{eqnarray}
where $\mathcal{Z}$ is an arbitrary PD tensor and $M_1 \geq 1$ is a positive constant. Given two PD random tensors $\mathcal{X} \in \mathbb{C}^{I_1 \times \cdots \times I_N \times I_1 \times \cdots \times I_N}$ and $\mathcal{Y} \in \mathbb{C}^{I_1 \times \cdots \times I_N \times I_1 \times \cdots \times I_N}$ and a determinsitic PD tensor $\mathcal{C} \in \mathbb{C}^{I_1 \times \cdots \times I_N \times I_1 \times \cdots \times I_N}$, if $\mathcal{X} \#_{x^{m}f(x)} \mathcal{Y}\preceq \mathcal{I}$ almost surely and $p \geq 1, q >0$, we have
\begin{eqnarray}\label{eq1:thm:4.1}
\mathrm{Pr}\left(\mathcal{X}^q \#_{x^{m}f(x)} \mathcal{Y}^q \npreceq \mathcal{C}\right) \leq \mathrm{Tr}\left(\mathbb{E}\left[\left(M_1\left(\prod_{k=1}^{m}\mathrm{K}_k\right)\mathcal{I}\right)^p\right]\star_N \mathcal{C}^{-1}\right),
\end{eqnarray}
where $\mathrm{K}_k$ are Kantorovich contants defined as
\begin{eqnarray}\label{eq1-1:thm:4.1}
\mathrm{K}_k \define \mathrm{K}(\lambda^{-1}_{\max}\left(\mathcal{X}^{-1}g^{(m-k)}(\mathcal{X})\right),\lambda^{-1}_{\min}\left(\mathcal{X}^{-1}g^{(m-k)}(\mathcal{X})\right),2q),
\end{eqnarray}
where $g(x)$ is defined by $g(x) \define \left(F^{[-1]}(x^{-1}) \right)^{-1}$, in which the function $F(x) \define x^{m-1}f(x)$.
\end{theorem}
\textbf{Proof:}
From Lemma~\ref{lma:odd 4.3} and Lemma~\ref{lma:even 4.4}, we have
\begin{eqnarray}
\mathcal{X}^q \#_{x^{m}f(x)} \mathcal{Y}^q \preceq M_1\left(\prod_{k=1}^{m}\mathrm{K}_k\right)\mathcal{I},
\end{eqnarray}
This thoerem is proved by Lemma~\ref{lma:Loewner ordering with Markov Cheb inequalities}
$\hfill\Box$

We have the following corollary from Theorem~\ref{thm:4.1} by replacing $f \in \mbox{TMI}^{1}$ to $f \in \mbox{TMD}^{1}$. 
\begin{corollary}\label{cor:thm:4.1}
Let $f \in \mbox{TMI}^{1}$ with pmd property and $m \in \mathbb{N}$ with $m \geq 2$. For  $q >0$, we also assume that 
\begin{eqnarray}\label{eq0:cor:thm:4.1}
M_2f(\mathcal{Z}^q)\succeq  f^q(\mathcal{Z}),
\end{eqnarray}
where $\mathcal{Z}$ is an arbitrary PD tensor and $M_2 \geq 1$ is a positive constant. Given two PD random tensors $\mathcal{X} \in \mathbb{C}^{I_1 \times \cdots \times I_N \times I_1 \times \cdots \times I_N}$ and $\mathcal{Y} \in \mathbb{C}^{I_1 \times \cdots \times I_N \times I_1 \times \cdots \times I_N}$ and a determinstic tensor $\mathcal{C} \in \mathbb{C}^{I_1 \times \cdots \times I_N \times I_1 \times \cdots \times I_N}$, if $\mathcal{X} \#_{x^{m}f(x)} \mathcal{Y}\succeq \mathcal{I}$ almost surely and $p \geq 1, q >0$, we have
\begin{eqnarray}\label{eq1:cor:thm:4.1}
\mathrm{Pr}\left(\left(M_2\prod_{k=1}^{m}\mathrm{K}_k\right)^{-1}\mathcal{I} \npreceq \mathcal{C}\right)
\leq \mathrm{Tr}\left(\mathbb{E}\left[\left(\mathcal{X}^q \#_{x^{m}f(x)} \mathcal{Y}^q\right)^p\right]\star_N \mathcal{C}^{-1}\right),
\end{eqnarray}
where $\mathrm{K}_k$ are Kantorovich contants defined as
\begin{eqnarray}\label{eq1-1:cor:thm:4.1}
\mathrm{K}_k \define \mathrm{K}(\lambda^{-1}_{\max}\left(\mathcal{X}^{-1}g^{(m-k)}(\mathcal{X})\right),\lambda^{-1}_{\min}\left(\mathcal{X}^{-1}g^{(m-k)}(\mathcal{X})\right),2q),
\end{eqnarray}
where $g(x)$ is defined by $g(x) \define \left(F^{[-1]}(x^{-1}) \right)^{-1}$, in which the function $F(x) \define x^{m-1}f(x)$.
\end{corollary}
\textbf{Proof:}
By setting $x^n f(x)$ as $x^n f^{-1}(1/x)$ in Theorem~\ref{thm:4.1}, we have
\begin{eqnarray}\label{eq1:cor:thm:4.1}
\mathcal{X}^q \#_{x^{m}f(x)} \mathcal{Y}^q \succeq \left(M_2\prod_{k=1}^{m}\mathrm{K}_k\right)^{-1}\mathcal{I}.
\end{eqnarray}
This thoerem is proved by Lemma~\ref{lma:Loewner ordering with Markov Cheb inequalities}
$\hfill\Box$

\section{Lie-Trotter Formula for Tensors and Its Applications}\label{sec:Lie-Trotter Formula and Its Applications} 

\subsection{Lie-Trotter Formula for Tensors}\label{sec:Lie-Trotter Formula for Tensors}

In this section, we will present a general Lie-Trotter formula for tensors with connection functions which are first differentiable functions on $(0, \infty)$. We use $\mathfrak{H}$ to represent the set of Hermitian tensors. We begin with the following Lemma.
\begin{lemma}\label{lma:5.2}
Let $g$ be a first differentialable function on $(0, \infty)$, a Hermitian tensor $\mathcal{X} \in \mathbb{C}^{I_1 \times \cdots \times I_N \times I_1 \times \cdots \times I_N}$, and a $\mathfrak{H}$-valued funtion $M(q)$ defined on $(-\epsilon_0, \epsilon_0)$, where $\epsilon >0$, such that $M(0)=0$ and $\lim\limits_{q\rightarrow 0}\frac{\left\Vert M(q)\right\Vert_{\rho}}{|q|}=0$, where $\left\Vert \cdot \right\Vert_{\rho}$ is an unitarily invariant norm defined by Eq.~\eqref{eq:def gauge func and general unitarily invariant norm}. Then, there exists another $\mathfrak{H}$-valued funtion $\tilde{M}(q)$ such that
\begin{eqnarray}
g(\mathcal{I}+q\mathcal{X}+M(q)) = f(1)\mathcal{I} + qg'(1)\mathcal{X} + \tilde{M}(q),
\end{eqnarray}
where $\tilde{M}(q)$ is defined on $(-\epsilon, \epsilon)$ for some $\epsilon \in (0, \epsilon_0)$; moreover,  
\begin{eqnarray}
\lim\limits_{q\rightarrow 0}\frac{\left\Vert \tilde{M}(q)\right\Vert_{\rho}}{|q|}=0.
\end{eqnarray}
\end{lemma}
\textbf{Proof:}
Because $\lim\limits_{q\rightarrow 0}\frac{\left\Vert M(q)\right\Vert_{\rho}}{|q|}=0$, we can find $a >0$ and 
$\epsilon \in (0, \epsilon_0)$ such that 
\begin{eqnarray}\label{eq1:lma:5.2}
\left\Vert \mathcal{X} + \frac{M(q)}{q}\right\Vert_{\rho} < a,
\end{eqnarray}
where $q \in (-\epsilon, \epsilon)\backslash \{0\}$. Then, for each $q \in (-\epsilon, \epsilon)\backslash \{0\}$, we can perform the spectral decomposition of $\mathcal{X} + \frac{M(q)}{q}$ as 
\begin{eqnarray}\label{eq2:lma:5.2}
\mathcal{X} + \frac{M(q)}{q} &=& \int_{-a'}^{a'}\lambda d E_q(\lambda),
\end{eqnarray}
where $[-a',a']$ is the eigenvalue range of the tensor $\mathcal{X} + \frac{M(q)}{q}$ given condition by Eq.~\eqref{eq1:lma:5.2}, and $E_q(\lambda)$ is eigen-tensor with respect to the eigenvalue $\lambda$ for the tensor $\mathcal{X} + \frac{M(q)}{q}$. From Eq.~\eqref{eq2:lma:5.2}, we can have the spectral decomposition for the tensor 
$g(\mathcal{I}+q\mathcal{X}+M(q))$ as 
\begin{eqnarray}\label{eq3:lma:5.2}
g(\mathcal{I}+q\mathcal{X}+M(q)) &=& \int_{-a'}^{a'}g(1+q \lambda) d E_q(\lambda). 
\end{eqnarray}

From the mean value therem, we have $g(1+q \lambda) = g(1) + q\lambda g'(1 +\theta q \lambda)$ for some $\theta \in (0,1)$. Then, we can express $g(1+q \lambda)$ as 
\begin{eqnarray}\label{eq4:lma:5.2}
g(1+q \lambda) &=& g(1) + qg'(1) + q\lambda \psi(q,\lambda),
\end{eqnarray}
where $\psi(q,\lambda) \define g'(1+\theta q \lambda) - g'(1)$. Note that $\sup\limits_{|\lambda|\leq a'}|\psi(q,\lambda)| \rightarrow 0$ as $q \rightarrow 0$. From Eq.~\eqref{eq3:lma:5.2} and Eq.~\eqref{eq4:lma:5.2}, we have
\begin{eqnarray}\label{eq5:lma:5.2}
g(\mathcal{I}+q\mathcal{X}+M(q))&=&g(1)\mathcal{I}+qg'(1)\left(\mathcal{X}+\frac{M(q)}{q}\right)+q \int_{-a'}^{a'}
\lambda\psi(q,\lambda)d E_q(\lambda).
\end{eqnarray}
From Eq.~\eqref{eq5:lma:5.2}, we have
\begin{eqnarray}\label{eq6:lma:5.2}
\left\Vert g(\mathcal{I}+q\mathcal{X}+M(q)) - g(1)\mathcal{I}-qg'(1)\mathcal{X}\right\Vert_{\rho}/|q|
&\leq& |g'(1)|\frac{\left\Vert M(q)\right\Vert_{\rho}}{|q|}+\sup\limits_{|\lambda|\leq a'}|\lambda\psi(q,\lambda)|.
\end{eqnarray}
Note that the R.H.S. of Eq.~\eqref{eq6:lma:5.2} will approach to zero as $q \rightarrow 0$ due to $\sup\limits_{|\lambda|\leq a'}|\psi(q,\lambda)| \rightarrow 0$ as $q \rightarrow 0$. Therefore, this Lemma is proved by setting  
$\tilde{M}(q)$ as
\begin{eqnarray}
\tilde{M}(q) \define g(\mathcal{I}+q\mathcal{X}+M(q)) - g(1)\mathcal{I}-qg'(1)\mathcal{X}.
\end{eqnarray}
$\hfill\Box$

We are ready to present the Lie-Trotter formula for tensors. 
\begin{theorem}\label{thm:5.1}
Let $g$ be a first differentiable function on $(0, \infty)$ with $g(1)=1$ and two Hermitian tensors $\mathcal{X}, \mathcal{Y} \in \mathbb{C}^{I_1 \times \cdots \times I_N \times I_1 \times \cdots \times I_N}$, then, we have
\begin{eqnarray}\label{eq1:thm:5.1}
\lim\limits_{q \rightarrow 0} \left(e^{q \mathcal{X}} \#_{g} e^{q \mathcal{Y}}\right)^{1/q} &=& \exp\left(g'(1)\mathcal{X}+(1-g'(1))\mathcal{Y}\right),
\end{eqnarray}
where $q \rightarrow 0$.
\end{theorem}
\textbf{Proof:}
From the Taylor expansions of $e^{q \mathcal{X}}$ and $e^{q \mathcal{Y}/2}$, we have 
\begin{eqnarray}
e^{-q \mathcal{Y}/2}\star_N e^{q \mathcal{X}} \star_N e^{-q \mathcal{Y}/2}
&=& \mathcal{I} + q(\mathcal{X} - \mathcal{Y}) + M(q),
\end{eqnarray}
where $M(q) \in \mathfrak{H}$ and $\lim\limits_{q \rightarrow 0} \frac{\left\Vert M(q)\right\Vert_{\rho}}{|q|}=0$. From Lemma~\ref{lma:5.2}, these exists a $\mathfrak{H}$-valued function $\left\Vert \tilde{M}(q)\right\Vert_{\rho}$ on $(-\epsilon, \epsilon)$ for some $\epsilon > 0$ such that 
\begin{eqnarray}
g(e^{-q \mathcal{Y}/2}\star_N e^{q \mathcal{X}} \star_N e^{-q \mathcal{Y}/2})
&=& \mathcal{I} + q g'(1)(\mathcal{X} - \mathcal{Y})+ \tilde{M}(q),
\end{eqnarray}
where $\lim\limits_{q \rightarrow 0} \frac{\left\Vert \tilde{M}(q)\right\Vert_{\rho}}{|q|}=0$. Then, we have
\begin{eqnarray}\label{eq2:thm:5.1}
e^{q \mathcal{X}}\#_{g}e^{q \mathcal{Y}}
&=& \left(\mathcal{I} + q g'(1)(\mathcal{X} - \mathcal{Y})+ \tilde{M}(q)\right)e^{q \mathcal{Y}},
\end{eqnarray}
by taking logarithm with respect to both sides of Eq.~\eqref{eq2:thm:5.1}, we have 
\begin{eqnarray}\label{eq3:thm:5.1}
\log(e^{q \mathcal{X}}\#_{g}e^{q \mathcal{Y}})
&=& q \mathcal{Y} + \log\left(\mathcal{I} + q g'(1)(\mathcal{X} - \mathcal{Y})+ \tilde{M}(q)\right) \nonumber \\
&=_1& q\left(g'(1)\mathcal{X}+(1-g'(1))\mathcal{Y}\right)+ \dbtilde{M}(q),
\end{eqnarray}
where we apply the approximation $\log(1+x) \approx x $ and note that the term $\dbtilde{M}(q)$ satsifies the following:  
\begin{eqnarray}\label{eq4:thm:5.1}
\lim\limits_{q \rightarrow 0}\frac{\left\Vert\dbtilde{M}(q)\right\Vert_{\rho}}{|q|} = 0.
\end{eqnarray}
Then, this theorem is proved by setting $q \rightarrow 0$ from Eq.~\eqref{eq3:thm:5.1}.
$\hfill\Box$

\subsection{Applications of Lie-Trotter Formula}\label{sec:Applications of Lie-Trotter Formula for Tensors}

In this section, we will apply the Lie-Trotter formula derived in Theorem~\ref{thm:5.1} to build new tail bounds of bivariate random tensor means for functions in $\mbox{TMI}^{1}$.

\begin{theorem}\label{thm:Lit-Trotter App TMI}
Let $0 < q \leq 1/2$, $g$ be a first differentialable function on $(0, \infty)$ with pmi property, $g(1)$ be equal to $1$, and $\mathcal{X}, \mathcal{Y} \in \mathbb{C}^{I_1 \times \cdots \times I_N \times I_1 \times \cdots \times I_N}$ be two Hermitian random tensors, then, we have the following tail bound:
\begin{eqnarray}\label{eq2-1:thm:Lit-Trotter App TMI}
\mathrm{Pr}\left(\exp\left((m+g'(1))\log\mathcal{X}+(1-m-g'(1))\log\mathcal{Y} \right) \npreceq \mathcal{C}\right)\leq
\mathrm{Tr}\left(\mathbb{E}\left[\left(\mathcal{X}^q \#_{x^{m}g(x)} \mathcal{Y}^q\right)^{r/q}\right]\star_N \mathcal{C}^{-1}\right),
\end{eqnarray}
where $r \geq 1$.

Let $0 < q \leq 1/2$, $f$ be a first differentialable function on $(0, \infty)$ with pmd property, $f(1)$ be equal to $1$, and $\mathcal{X}, \mathcal{Y} \in \mathbb{C}^{I_1 \times \cdots \times I_N \times I_1 \times \cdots \times I_N}$ be two Hermitian random tensors, then, we have the following tail bound:
\begin{eqnarray}\label{eq2-2:thm:Lit-Trotter App TMI}
\mathrm{Pr}\left(\left(\mathcal{X}^q \#_{x^{m}f(x)} \mathcal{Y}^q\right)^{1/q}\npreceq \mathcal{C}\right)\leq
\mathrm{Tr}\left(\mathbb{E}\left[\left((m+f'(1))\log\mathcal{X}+(1-m-f'(1))\log\mathcal{Y} \right)^{r}\right]\star_N \mathcal{C}^{-1}\right),
\end{eqnarray}
where $r \geq 1$.
\end{theorem}
\textbf{Proof:}
From the proof in Theorem~\ref{thm:5.1} and the condition $0 < q \leq 1/2$, we have
\begin{eqnarray}\label{eq3:thm:Lit-Trotter App TMI}
\mathcal{X}^q \#_{x^{m}g(x)} \mathcal{Y}^q \preceq \mathcal{I},
\end{eqnarray}
which is equivalent to the following 
\begin{eqnarray}\label{eq4:thm:Lit-Trotter App TMI}
\left(\mathcal{X}^p \#_{x^{m}g(x)} \mathcal{Y}^p\right)^{1/p} 
\preceq \left(\mathcal{X}^q \#_{x^{m}g(x)} \mathcal{Y}^q\right)^{1/q},
\end{eqnarray}
where $0 < p \leq q/2$. By taking $p \rightarrow 0$ and applying Theorem~\ref{thm:5.1} in Eq.~\eqref{eq4:thm:Lit-Trotter App TMI}, we have
\begin{eqnarray}\label{eq5:thm:Lit-Trotter App TMI}
\exp\left((m+g'(1))\log\mathcal{X}+(1-m-g'(1))\log\mathcal{Y} \right)\preceq  \left(\mathcal{X}^q \#_{x^{m}g(x)} \mathcal{Y}^q\right)^{1/q}.
\end{eqnarray}
Applying Lemma~\ref{lma:Loewner ordering with Markov Cheb inequalities} to Eq.~\eqref{eq5:thm:Lit-Trotter App TMI}, we have Eq.~\eqref{eq2-1:thm:Lit-Trotter App TMI}.

From the proof in Corollary~\ref{cor:thm:4.1} and the condition $0 < q \leq 1/2$, we have
\begin{eqnarray}\label{eq6:thm:Lit-Trotter App TMI}
\mathcal{X}^q \#_{x^{m}f(x)} \mathcal{Y}^q \succeq \mathcal{I},
\end{eqnarray}
which is equivalent to the following 
\begin{eqnarray}\label{eq7:thm:Lit-Trotter App TMI}
\left(\mathcal{X}^p \#_{x^{m}f(x)} \mathcal{Y}^p\right)^{1/p} 
\succeq \left(\mathcal{X}^q \#_{x^{m}f(x)} \mathcal{Y}^q\right)^{1/q},
\end{eqnarray}
where $0 < p \leq q/2$. By taking $p \rightarrow 0$ and applying Theorem~\ref{thm:5.1} in Eq.~\eqref{eq7:thm:Lit-Trotter App TMI}, we have
\begin{eqnarray}\label{eq8:thm:Lit-Trotter App TMI}
\exp\left((m+f'(1))\log\mathcal{X}+(1-m-f'(1))\log\mathcal{Y} \right)\succeq  \left(\mathcal{X}^q \#_{x^{m}f(x)} \mathcal{Y}^q\right)^{1/q}.
\end{eqnarray}
Applying Lemma~\ref{lma:Loewner ordering with Markov Cheb inequalities} to Eq.~\eqref{eq5:thm:Lit-Trotter App TMI}, we have Eq.~\eqref{eq2-2:thm:Lit-Trotter App TMI}.
$\hfill\Box$

\section{L\"oewner and Majorization Ordering for Bivariate Random Tensor Means with Non-invertible Tensors}\label{sec:Loewner and Majorization Ordering for Bivariate Random Tensor Means with Non-invertible Tensors}

\subsection{Limitation Method by Random Tensor Topology (RTT)}\label{sec:Limitation Method by Random Tensor Topology (RTT)}

Previous discussions until now are based on tensors with PD property. In this section, we will apply the \emph{random tensor topology} (RTT) concept to deal with PSD tensors from PD tensors by taking the limit operation. The notion about RTT is provided by the following definition.
\begin{definition}\label{def:conv in mean}
We say that a sequence of random tensor $\mathcal{X}_n$ converges to the random tensor $\mathcal{X}$ with respect to the tensor norm $\left\Vert \cdot \right\Vert_{\rho}$ in the sense of RTT, if we have
\begin{eqnarray}
\mathbb{E}\left( \left\Vert \mathcal{X}_n \right\Vert_{\rho} \right)~~~\mbox{exists,}
\end{eqnarray}
and
\begin{eqnarray}
\lim\limits_{n \rightarrow \infty}\mathbb{E}\left( \left\Vert \mathcal{X}_n  - \mathcal{X} \right\Vert_{\rho} \right) = 0. 
\end{eqnarray}
We adopt the notation $\lim\limits_{n \rightarrow \infty}\mathcal{X}_n = \mathcal{X}$ to represent that random tensors $\mathcal{X}_n$ converges to the random tensor $\mathcal{X}$ with respect to the tensor norm $\left\Vert \cdot \right\Vert_{\rho}$ in the sense of RTT.
\end{definition}

Given two PSD tensors $\mathcal{X}, \mathcal{Y}$, we will determine the term $\mathcal{X} \#_{g} \mathcal{Y}$ by $\lim\limits_{\epsilon \rightarrow 0^+} \left(\mathcal{X}+\epsilon \mathcal{I}\right)\#_{g} \left(\mathcal{Y}+\epsilon \mathcal{I}\right)$ under RTT. We will begin with several lemmas which will be used for later theorem proof.

\begin{lemma}\label{lma:prop 2.2}
For any function $g > 0$, set $h(x)\define xg(x^{-1})$, we have the following equivalent conditions:
\begin{enumerate}[label=(\roman*)]
\item $g \in \mbox{TC}$ and $g(0^{+})=0$;
\item $g(0^{+})=0$ and the connection function $g$ is jointly tensor convex, i.e., 
\begin{eqnarray}\label{eq1:lma:prop 2.2}
\left(\lambda\mathcal{X}_1+\overline{\lambda}\mathcal{X}_2\right)\#_g \left(\lambda\mathcal{Y}_1+\overline{\lambda}\mathcal{Y}_2\right)\preceq \lambda\mathcal{X}_1 \#_g \mathcal{Y}_1+\overline{\lambda}\mathcal{X}_2 \#_g \mathcal{Y}_2,
\end{eqnarray}
where $\overline{\lambda}\define 1 - \lambda$, and $\mathcal{X}_1,\mathcal{X}_2,\mathcal{Y}_1,\mathcal{Y}_2 \in \mathbb{C}^{I_1 \times \cdots \times I_N \times I_1 \times \cdots \times I_N}$ are PD tensors;
\item Given any PD tensor $\mathcal{X} \in \mathbb{C}^{I_1 \times \cdots \times I_N \times I_1 \times \cdots \times I_N}$, the operation $\#_g$ is right tensor decreasing, i.e., 
\begin{eqnarray}\label{eq2:lma:prop 2.2}
\mathcal{O} \prec \mathcal{Y}_1 \preceq \mathcal{Y}_2 \Rightarrow \mathcal{X} \#_{g} \mathcal{Y}_1 \succeq 
\mathcal{X} \#_{g} \mathcal{Y}_2;
\end{eqnarray}
\item Given any PD tensor $\mathcal{Y} \in \mathbb{C}^{I_1 \times \cdots \times I_N \times I_1 \times \cdots \times I_N}$, the operation $\#_h$ is left tensor decreasing, i.e., 
\begin{eqnarray}\label{eq3:lma:prop 2.2}
\mathcal{O} \prec \mathcal{X}_1 \preceq \mathcal{X}_2 \Rightarrow \mathcal{X}_1 \#_{h} \mathcal{Y} \succeq 
\mathcal{X}_2 \#_{h} \mathcal{Y}.
\end{eqnarray}
\end{enumerate}
\end{lemma}
\textbf{Proof:}
The equivalence between (i) and (ii) can be found at Theorem 2.2 in~\cite{ebadian2011perspectives} by treating tensors as operators. Because $h\define xg(x^{-1})$, we have the equivalence between (i) and (iv) due to 
\begin{eqnarray}
h(\mathcal{X})=\mathcal{Y}^{-1/2}\star_N \left((\mathcal{Y}^{1/2}\star_N \mathcal{X} \star_N \mathcal{Y}^{1/2})\#_h \mathcal{Y}\right) \star_N \mathcal{Y}^{-1/2}.
\end{eqnarray}
Finally, the equivalence between (iii) and (iv) comes from the following relationship since $h(x) \define xg(x^{-1})$:
\begin{eqnarray}
\mathcal{X} \#_{h} \mathcal{Y} &=& \mathcal{Y} \#_{g} \mathcal{X}.
\end{eqnarray}
$\hfill \Box$

Given two PSD tensors $\mathcal{X}, \mathcal{Y}$ with $\mathcal{X} \preceq c \mathcal{Y}$ for some $c>0$, there is a unique PD tensor, denoted by $\eta(\mathcal{X}, \mathcal{Y}) \in \mathbb{C}^{I_1 \times \cdots \times I_N \times I_1 \times \cdots \times I_N}$, such that $\eta(\mathcal{X}, \mathcal{Y})(\mathcal{I} - \mbox{PR}(\mathcal{Y})) = \mathcal{O}$ and $\mathcal{X} = \mathcal{Y}^{1/2}\star_N \eta(\mathcal{X}, \mathcal{Y}) \star_N \mathcal{Y}^{1/2}$, where $\mbox{PR}(\mathcal{Y})$ is the orthogonal projection onto the the closure of the range of $\mathcal{Y}$.  Note that $\eta(\mathcal{X}, \mathcal{Y})$ will be reduced as $\mathcal{Y}^{-1/2}\star_N \mathcal{X} \star_N \mathcal{Y}^{-1/2}$ when $\mathcal{Y} \succ \mathcal{O}$. The next lemma will show the existence of $\lim\limits_{\epsilon \rightarrow 0^+} \left(\mathcal{X}+\epsilon \mathcal{I}\right)\#_{g} \left(\mathcal{Y}+\epsilon \mathcal{I}\right)$.
\begin{lemma}\label{lma:thm 6.2}
Let $g \in \mbox{TC}$, then we have the following equivalent conditions:
\begin{enumerate}[label=(\roman*)]
\item For every PSD random tensors $\mathcal{X}, \mathcal{Y}$ with $\mathcal{X} \preceq c \mathcal{Y}$ for some $c >0$, we have
\begin{eqnarray}
\lim\limits_{\epsilon \rightarrow 0^+} \left(\mathcal{X}+\epsilon\mathcal{I}\right)\#_{g} \left(\mathcal{Y}+\epsilon \mathcal{I}\right)&=&\mathcal{Y}^{1/2}\star_N g(\eta(\mathcal{X},\mathcal{Y})) \star_N \mathcal{Y}^{1/2};
\end{eqnarray}
\item $g(0^{+}) < \infty$.
\end{enumerate}
\end{lemma}
\textbf{Proof:}
For the direction, (i) $\Rightarrow$ (ii), by setting $\mathcal{X}=x \mathcal{I}$ and $\mathcal{Y}=y \mathcal{I}$ with $x,y \geq 0$, we have
\begin{eqnarray}
 \left(\mathcal{X}+\epsilon\mathcal{I}\right)\#_{g} \left(\mathcal{Y}+\epsilon \mathcal{I}\right)
&=& (y+\epsilon)g\left(\frac{x+\epsilon}{y+\epsilon}\right)\mathcal{I}.
\end{eqnarray}
If we set $x=0$ and $y=1$, we have $(1+\epsilon)g(\frac{\epsilon}{1+\epsilon}) \rightarrow g(0^+)$ as $\epsilon \rightarrow 0^+$. 

We now prove the direction, (ii) $\Rightarrow$ (i). From Theorem 8.1 in~\cite{hiai2011quantum}, we can represent a tensor (operator) convex function $g(x)$ as
\begin{eqnarray}
g(x)&=&a_0+a_1x+a_2x^2 + \int_{0}^{\infty}\left(\frac{x}{1+s} - \frac{x}{x+s}\right)d \mu(s),
\end{eqnarray} 
where $a_0, a_1 \in \mathbb{R}$, $a_2 >0$, and $\mu(s)$ is a positive measure on $(0, \infty)$ satisfying 
\begin{eqnarray}
\int_{0}^{\infty}(1+s)^{-2} < \infty.
\end{eqnarray}

By setting
\begin{eqnarray}
\psi_s(x) \define \left(\frac{x}{1+s} - \frac{x}{x+s}\right), 
\end{eqnarray}
we have
\begin{eqnarray}\label{eq6-4:lma:thm 6.2}
\left(\mathcal{X}+\epsilon\mathcal{I}\right)\#_{g} \left(\mathcal{Y}+\epsilon \mathcal{I}\right)
&=& a_0 \left(\mathcal{Y}+\epsilon \mathcal{I}\right)+ a_1 \left(\mathcal{X}+\epsilon \mathcal{I}\right)
+a_2 \left(\mathcal{X}+\epsilon\mathcal{I}\right)\#_{x^2} \left(\mathcal{Y}+\epsilon \mathcal{I}\right) \nonumber \\
&&+ \int_{0}^{\infty}\left(\mathcal{X}+\epsilon\mathcal{I}\right)\#_{\psi_s(x)} \left(\mathcal{Y}+\epsilon \mathcal{I}\right)d \mu(s).
\end{eqnarray}

From the assumption that $\mathcal{X} \preceq c\mathcal{Y}$ for any $c >0$,  we have $\frac{cx+\epsilon}{x+\epsilon} \leq \tilde{c}$ for all $x \geq 0$ and any $\epsilon >0$, where $\tilde{c}=c$ if $c \geq 1$ and $\tilde{c}=1$ if $c < 1$. Then, we have
\begin{eqnarray}\label{eq1:dom bound lma:thm 6.2}
\left(\mathcal{Y}+\epsilon \mathcal{I}\right)^{-1/2}\left(\mathcal{X}+\epsilon \mathcal{I}\right)\left(\mathcal{Y}+\epsilon \mathcal{I}\right)^{-1/2}
&\preceq&\left(\mathcal{Y}+\epsilon \mathcal{I}\right)^{-1/2}\left(c\mathcal{Y}+\epsilon \mathcal{I}\right)\left(\mathcal{Y}+\epsilon \mathcal{I}\right)^{-1/2} \nonumber \\
&=& \left(c\mathcal{Y}+\epsilon \mathcal{I}\right)\left(\mathcal{Y}+\epsilon \mathcal{I}\right)^{-1}  \nonumber \\
&\preceq& \tilde{c} \mathcal{I},
\end{eqnarray}
Besides, we also have the following bounds for the function $\psi_s(x)$ as
\begin{eqnarray}
\frac{-1}{(1+s)^2} \leq \psi_s(x) \leq \psi_s(\tilde{c}),
\end{eqnarray}
where $x \in [0, \tilde{c}]$. Then, we have
\begin{eqnarray}\label{eq2:dom bound lma:thm 6.2}
\frac{-1}{(1+s)^2}\mathcal{I} \preceq \psi_s(\left(\mathcal{Y}+\epsilon \mathcal{I}\right)^{-1/2}\left(\mathcal{X}+\epsilon \mathcal{I}\right)\left(\mathcal{Y}+\epsilon \mathcal{I}\right)^{-1/2}) \preceq \psi_s(\tilde{c})\mathcal{I}.
\end{eqnarray}
If $\epsilon \in (0,1)$, Eq.~\eqref{eq2:dom bound lma:thm 6.2} implies the following bound 
\begin{eqnarray}\label{eq3:dom bound lma:thm 6.2}
\frac{-\left(\lambda_{\max}(\mathcal{Y})+1\right)}{(1+s)^2}\mathcal{I} \preceq \psi_s(\left(\mathcal{Y}+\epsilon \mathcal{I}\right)^{-1/2}\left(\mathcal{X}+\epsilon \mathcal{I}\right)\left(\mathcal{Y}+\epsilon \mathcal{I}\right)^{-1/2}) \preceq \psi_s(\tilde{c})\left(\lambda_{\max}(\mathcal{Y})+1\right)\mathcal{I}.
\end{eqnarray}

Suppose that the following two limits exist in the sense of RTT:
\begin{eqnarray}\label{eq6-5:lma:thm 6.2}
\lim\limits_{\epsilon \rightarrow 0^+}(\mathcal{X}+\epsilon \mathcal{I})\#_{x^2}(\mathcal{Y}+\epsilon \mathcal{I}) &=& \mathcal{Y}^{1/2}\eta^2(\mathcal{X}, \mathcal{Y})\mathcal{Y}^{1/2},
\end{eqnarray}
\begin{eqnarray}\label{eq6-6:lma:thm 6.2}
\lim\limits_{\epsilon \rightarrow 0^+}(\mathcal{X}+\epsilon \mathcal{I})\#_{ \psi_s(x)}(\mathcal{Y}+\epsilon \mathcal{I}) &=& \mathcal{Y}^{1/2}\psi_s(\eta(\mathcal{X}, \mathcal{Y}))\mathcal{Y}^{1/2};
\end{eqnarray}
then, from Eq.~\eqref{eq3:dom bound lma:thm 6.2}, $\int_0^{\infty}(1+s)^{-2} d \mu(s) < \infty$, $\int_0^{\infty}\psi_s (\tilde{c}) d \mu(s) < \infty$, and Lebesgue convergence theorem, we have 
\begin{eqnarray}\label{eq6-7:lma:thm 6.2}
\int_0^{\infty}\mathcal{Y}^{1/2}\psi_s(\eta(\mathcal{X}, \mathcal{Y}))\mathcal{Y}^{1/2} d \mu(s) = \lim\limits_{\epsilon \rightarrow 0^+}\int_0^{\infty}(\mathcal{X}+\epsilon \mathcal{I})\#_{ \psi_s(x)}(\mathcal{Y}+\epsilon \mathcal{I}) d \mu(s).
\end{eqnarray}
Therefore, from Eqs.~\eqref{eq6-4:lma:thm 6.2},~\eqref{eq6-5:lma:thm 6.2} and~\eqref{eq6-7:lma:thm 6.2}, we prove the statement (i) by:
\begin{eqnarray}\label{eq6-8:lma:thm 6.2}
\lim\limits_{\epsilon \rightarrow 0^+}\left(\mathcal{X}+\epsilon\mathcal{I}\right)\#_{g} \left(\mathcal{Y}+\epsilon \mathcal{I}\right)
&=& a_0 \mathcal{Y}+ a_1 \mathcal{X}+a_2 \mathcal{Y}^{1/2}\eta^2(\mathcal{X}, \mathcal{Y})\mathcal{Y}^{1/2}\nonumber \\
&  & + \int_{0}^{\infty}\mathcal{Y}^{1/2}\psi_s(\eta(\mathcal{X}, \mathcal{Y}))\mathcal{Y}^{1/2}d \mu(s) \nonumber \\
&=&\mathcal{Y}^{1/2}\left(a_0 \mathcal{I}+a_1 \eta(\mathcal{X}, \mathcal{Y}) + a_2 \eta^2(\mathcal{X}, \mathcal{Y})\right. \nonumber \\
&   &\left.+ \int_{0}^{\infty}\psi_s(\eta(\mathcal{X}, \mathcal{Y}))d \mu(s) \right)\mathcal{Y}^{1/2}\nonumber \\
&=&\mathcal{Y}^{1/2}g(\eta(\mathcal{X},\mathcal{Y}))\mathcal{Y}^{1/2}.
\end{eqnarray}

Our final tasks are to prove Eq.~\eqref{eq6-5:lma:thm 6.2} and Eq.~\eqref{eq6-6:lma:thm 6.2}. Because $\mathcal{X}\preceq c \mathcal{Y}$, we can have a bounded tensor $\mathcal{U}$ with $\left\Vert\mathcal{U}\right\Vert_{\rho} \leq c^{1/2}$ such that $\mathcal{U}(\mathcal{I}-\mbox{PR}(\mathcal{Y}))=\mathcal{O}$ and $\mathcal{X}^{1/2}=\mathcal{U}\mathcal{Y}^{1/2}=\mathcal{Y}^{1/2}\mathcal{U}^{\mathrm{H}}$. Then, we have $\eta(\mathcal{X}, \mathcal{Y})=\mathcal{U}^{\mathrm{H}}\mathcal{U}$ and we can express the tensor $\mathcal{X}$ by  $\eta(\mathcal{X}, \mathcal{Y})$ as
\begin{eqnarray}\label{eq6-8-1:lma:thm 6.2}
\mathcal{X}&=&\mathcal{Y}^{1/2}\eta(\mathcal{X}, \mathcal{Y})\mathcal{Y}^{1/2}.
\end{eqnarray}

To prove Eq.~\eqref{eq6-5:lma:thm 6.2}, we will expand $(\mathcal{X}+\epsilon \mathcal{I})\#_{x^2}(\mathcal{Y}+\epsilon \mathcal{I})$ as follows:
\begin{eqnarray}\label{eq6-8-2:lma:thm 6.2}
(\mathcal{X}+\epsilon \mathcal{I})\#_{x^2}(\mathcal{Y}+\epsilon \mathcal{I})&=&
(\mathcal{X}+\epsilon \mathcal{I})(\mathcal{Y}+\epsilon \mathcal{I})^{-1}(\mathcal{X}+\epsilon \mathcal{I})\nonumber \\
&=& \underbrace{\mathcal{X}(\mathcal{Y}+\epsilon \mathcal{I})^{-1}\mathcal{X}}_A+\underbrace{\epsilon\mathcal{X}(\mathcal{Y}+\epsilon \mathcal{I})^{-1}}_B\nonumber \\
&  &+\underbrace{\epsilon(\mathcal{Y}+\epsilon \mathcal{I})^{-1}\mathcal{X}}_C+ 
\underbrace{\epsilon^2(\mathcal{Y}+\epsilon \mathcal{I})^{-1}\mathcal{I}}_D.
\end{eqnarray}
For the term corresponding to $B$, if we have the spectral decomposition of $\mathcal{Y}=\int_0^{\lambda_{\max}(\mathcal{Y})}\lambda d E(\lambda)$, we have 
\begin{eqnarray}\label{eq6-8-3:lma:thm 6.2}
\mathbb{E}\left[\left\Vert\epsilon\mathcal{X}(\mathcal{Y}+\epsilon \mathcal{I})^{-1}\right\Vert_{\rho}\right]
&=_1&\mathbb{E}\left[\left\Vert\epsilon\mathcal{Y}^{1/2}\eta(\mathcal{X}, \mathcal{Y})\mathcal{Y}^{1/2}(\mathcal{Y}+\epsilon \mathcal{I})^{-1}\right\Vert_{\rho}\right]\nonumber \\
&\leq&c\mathbb{E}\left[\lambda_{\max}(\mathcal{Y}^{1/2})\left\Vert\epsilon\mathcal{Y}^{1/2}(\mathcal{Y}+\epsilon \mathcal{I})^{-1}\right\Vert_{\rho}\right] \nonumber \\
&\leq&c\mathbb{E}\left[\lambda_{\max}(\mathcal{Y}^{1/2})\int_0^{\lambda_{\max}(\mathcal{Y})}\frac{\epsilon\lambda^{1/2}}{\lambda+\epsilon}
d\left\Vert E(\lambda)\right\Vert_{\rho}\right],
\end{eqnarray} 
where $=_1$ comes from Eq.~\eqref{eq6-8-1:lma:thm 6.2}. From Eq.~\eqref{eq6-8-3:lma:thm 6.2}, since $\frac{\epsilon\lambda^{1/2}}{\lambda+\epsilon} \rightarrow 0$ as $\epsilon \rightarrow 0^+$, we have $\left\Vert\epsilon\mathcal{X}(\mathcal{Y}+\epsilon \mathcal{I})^{-1}\right\Vert_{\rho} \rightarrow 0$ as $\epsilon \rightarrow 0^+$, i.e., $\epsilon\mathcal{X}(\mathcal{Y}+\epsilon \mathcal{I})^{-1} \rightarrow \mathcal{O}$ in RTT as $\epsilon \rightarrow 0^+$. Similar arguments used by Eq.~\eqref{eq6-8-3:lma:thm 6.2} can show that the terms $C$ and $D$ in Eq.~\eqref{eq6-8-2:lma:thm 6.2} approach to zero tensor in RTT as $\epsilon \rightarrow 0^+$. 

For the term corresponding to $A$, we have
\begin{eqnarray}\label{eq6-8-4:lma:thm 6.2}
\lim\limits_{\epsilon \rightarrow 0^+}\mathcal{X}(\mathcal{Y}+\epsilon \mathcal{I})^{-1}\mathcal{X}
&=& \lim\limits_{\epsilon \rightarrow 0^+}\mathcal{X}^{1/2}\mathcal{U}\mathcal{Y}^{1/2}(\mathcal{Y}+\epsilon \mathcal{I})^{-1}\mathcal{Y}^{1/2}\mathcal{U}^{\mathrm{H}}\mathcal{X}^{1/2}\nonumber \\
&=& \lim\limits_{\epsilon \rightarrow 0^+}\mathcal{X}^{1/2}\mathcal{U}\mathcal{Y}(\mathcal{Y}+\epsilon \mathcal{I})^{-1}\mathcal{U}^{\mathrm{H}}\mathcal{X}^{1/2}\nonumber \\
&=_1&\mathcal{X}^{1/2}\mathcal{U}\mbox{PR}(\mathcal{Y})\mathcal{U}^{\mathrm{H}}\mathcal{X}^{1/2}\nonumber \\
&=& \mathcal{Y}^{1/2}\mathcal{U}^{\mathrm{H}}\mathcal{U}\mathcal{U}^{\mathrm{H}}\mathcal{U}\mathcal{Y}^{1/2}\nonumber \\
&=&\mathcal{Y}^{1/2}\eta^2(\mathcal{X}, \mathcal{Y})\mathcal{Y}^{1/2}
\end{eqnarray} 
where $=_1$ comes from $\mathcal{Y}(\mathcal{Y}+\epsilon \mathcal{I})^{-1} \rightarrow \mbox{PR}(\mathcal{Y})$  in RTT as $\epsilon \rightarrow 0^+$.
Hence, we prove Eq.~\eqref{eq6-5:lma:thm 6.2}.

To prove Eq.~\eqref{eq6-6:lma:thm 6.2}, we have the following expanstion for $(\mathcal{X}+\epsilon \mathcal{I})\#_{ \psi_s(x)}(\mathcal{Y}+\epsilon \mathcal{I})$, which is 
\begin{eqnarray}\label{eq6-9-1:lma:thm 6.2}
(\mathcal{X}+\epsilon \mathcal{I})\#_{ \psi_s(x)}(\mathcal{Y}+\epsilon \mathcal{I})
&=&\underbrace{\frac{\mathcal{X}+\epsilon\mathcal{I}}{1+s}}_A- \nonumber \\
&   &\underbrace{\left[(\mathcal{Y}+\epsilon\mathcal{I})^{-1/2}(\mathcal{X}+\epsilon\mathcal{I})(\mathcal{Y}+\epsilon\mathcal{I})^{-1/2}+ s\mathcal{I}\right]^{-1}(\mathcal{X}+\epsilon\mathcal{I})}_B,
\end{eqnarray}
and, by taking $\epsilon \rightarrow 0$ at Eq~\eqref{eq6-9-1:lma:thm 6.2}, we claim Eq~\eqref{eq6-9-1:lma:thm 6.2} becomes:
\begin{eqnarray}\label{eq6-9-2:lma:thm 6.2}
(\mathcal{X}+\epsilon \mathcal{I})\#_{ \psi_s(x)}(\mathcal{Y}+\epsilon \mathcal{I})
&=& \mathcal{Y}^{1/2}\frac{\eta(\mathcal{X}, \mathcal{Y})}{1+s}\mathcal{Y}^{1/2}
-\mathcal{Y}^{1/2}(\eta(\mathcal{X}, \mathcal{Y})+s\mathcal{I})^{-1}\eta(\mathcal{X}, \mathcal{Y})\mathcal{Y}^{1/2} \nonumber \\
&=& \mathcal{Y}^{1/2}\psi_s(\eta(\mathcal{X}, \mathcal{Y}))\mathcal{Y}^{1/2}.
\end{eqnarray}
The limitation of the part A in Eq.~\eqref{eq6-9-1:lma:thm 6.2} becomes $\mathcal{Y}^{1/2}\frac{\eta(\mathcal{X}, \mathcal{Y})}{1+s}\mathcal{Y}^{1/2}$ can be obtained from Eq.~\eqref{eq6-8-1:lma:thm 6.2}. For the part B in Eq.~\eqref{eq6-9-1:lma:thm 6.2}, we consider 
\begin{eqnarray}\label{eq6-9-3:lma:thm 6.2}
\lim\limits_{\epsilon \rightarrow 0}\mathbb{E}\left[\left\Vert\left((\mathcal{Y}+\epsilon\mathcal{I})^{-1/2}(\mathcal{X}+\epsilon\mathcal{I})(\mathcal{Y}+\epsilon\mathcal{I})^{-1/2}+ s\mathcal{I}\right)^{-1}(\mathcal{X}+\epsilon\mathcal{I}) - (\eta(\mathcal{X}, \mathcal{Y})+ s\mathcal{I})^{-1}\mathcal{X}\right\Vert_{\rho}\right]\nonumber \\
=_1\lim\limits_{\epsilon \rightarrow 0}\mathbb{E}\left[\left\Vert\left(\eta(\mathcal{X}+\epsilon\mathcal{I}, \mathcal{Y}+\epsilon\mathcal{I})+ s\mathcal{I}\right)^{-1}(\mathcal{X}+\epsilon\mathcal{I}) - (\eta(\mathcal{X}, \mathcal{Y})+ s\mathcal{I})^{-1}\mathcal{X}\right\Vert_{\rho}\right]~~~~~~~~~~~~~~~~~~~~~~~~~~ \nonumber \\
\leq_2\lim\limits_{\epsilon \rightarrow 0}\mathbb{E}\left[\frac{\epsilon\left\Vert\eta(\mathcal{X}, \mathcal{Y})\right\Vert_{\rho}+
\epsilon s \left\Vert \mathcal{I}\right\Vert_{\rho}+\left\Vert\mathcal{X}\right\Vert_{\rho}\left\Vert \eta(\mathcal{X}+\epsilon\mathcal{I}, \mathcal{Y}+\epsilon\mathcal{I}) - \eta(\mathcal{X},\mathcal{Y})\right\Vert_{\rho}}{
(s+\lambda_{\min}(\eta(\mathcal{X}+\epsilon\mathcal{I}, \mathcal{Y}+\epsilon\mathcal{I})))(s+\lambda_{\min}(\eta(\mathcal{X}, \mathcal{Y})))\left\Vert \mathcal{I} \right\Vert_{\rho}}\right] 
=_3 0  ~~~~~~~~~~~~~~~
\end{eqnarray}
where $=_1$ comes from Eq.~\eqref{eq6-8-1:lma:thm 6.2}, $\leq_2$ comes from triangle inequality of the norm $\left\Vert \cdot \right\Vert_{\rho}$ and the lower bound of the norm $\left(\eta(\mathcal{X}+\epsilon\mathcal{I}, \mathcal{Y}+\epsilon\mathcal{I})+ s\mathcal{I}\right)(\eta(\mathcal{X}, \mathcal{Y})+ s\mathcal{I})$, $=_3$ comes from that each term in the numerator is bounded. From Eq.~\eqref{eq6-9-3:lma:thm 6.2}, the part B in Eq.~\eqref{eq6-9-1:lma:thm 6.2} will be $(\mathcal{Y}^{-1/2}\mathcal{X}\mathcal{Y}^{-1/2}+ s\mathcal{I})^{-1}\mathcal{X}$ by taking $\epsilon \rightarrow 0$. Therefore, we prove Eq.~\eqref{eq6-6:lma:thm 6.2}.
$\hfill \Box$

We are ready to present the following theorem about the limitation of the bivariate tensor mean by using a sequence of PD tensors. 
\begin{theorem}\label{thm:6.3}
Given a function $g \in \mbox{TC}$ defined on $(0, \infty)$ with $g(0^+) < \infty$, and two PSD random tensors $\mathcal{X}, \mathcal{Y} \in \mathbb{C}^{I_1 \times \dots \times I_N \times I_1 \times \dots \times I_N}$ with $\mathcal{X} \preceq c \mathcal{Y}$ almost surely for some $c >0$. If a sequence of PD tensors $\{\mathcal{A}_n\}$ satisfy $\left\Vert \mathcal{A}_n \right\Vert_{\rho} \rightarrow 0$, then, we have
\begin{eqnarray}\label{eq1:thm:6.3}
\lim\limits_{n \rightarrow \infty} \mathcal{X} \#_{g} (\mathcal{Y}+\mathcal{A}_n) &=&  
\lim\limits_{\epsilon \rightarrow 0^+} \mathcal{X} \#_{g} (\mathcal{Y}+\epsilon \mathcal{I}) \nonumber \\
&=&\mathcal{Y}^{1/2}\star_N g(\eta(\mathcal{X}, \mathcal{Y}))\star_N \mathcal{Y}^{1/2}.
\end{eqnarray}
\end{theorem}
\textbf{Proof:}
From Lemma~\ref{lma:thm 6.2}, we have 
\begin{eqnarray}\label{eq2:thm:6.3}
\lim\limits_{\epsilon \rightarrow 0^+} \mathcal{X} \#_{g} (\mathcal{Y}+\epsilon \mathcal{I}) &=&
\mathcal{Y}^{1/2}g(\eta(\mathcal{X}, \mathcal{Y}))\mathcal{Y}^{1/2}.
\end{eqnarray}

For the term $\lim\limits_{n \rightarrow \infty} \mathcal{X} \#_{g} (\mathcal{Y}+\mathcal{A}_n)$, we have
\begin{eqnarray}\label{eq3:thm:6.3}
\lim\limits_{n \rightarrow \infty} \mathcal{X} \#_{g} (\mathcal{Y}+\mathcal{A}_n) &=&
\lim\limits_{n \rightarrow \infty} (\mathcal{Y}+\mathcal{A}_n)^{1/2}g((\mathcal{Y}+\mathcal{A}_n)^{-1/2}\mathcal{X} (\mathcal{Y}+\mathcal{A}_n)^{-1/2} )(\mathcal{Y}+\mathcal{A}_n)^{1/2} \nonumber \\
&=_1&\lim\limits_{n \rightarrow \infty} (\mathcal{Y}+\mathcal{A}_n)^{1/2}g((\mathcal{Y}+\mathcal{A}_n)^{-1/2}\mathcal{Y}^{1/2}\eta(\mathcal{X}, \mathcal{Y})\mathcal{Y}^{1/2}(\mathcal{Y}+\mathcal{A}_n)^{-1/2})\nonumber \\
&& \star_N (\mathcal{Y}+\mathcal{A}_n)^{1/2} \nonumber \\
&=_2&\mathcal{Y}^{1/2}g(\mbox{PR}(\mathcal{Y})\eta(\mathcal{X}, \mathcal{Y})\mbox{PR}(\mathcal{Y}))\mathcal{Y}^{1/2}=\mathcal{Y}^{1/2}g(\eta(\mathcal{X}, \mathcal{Y}))\mathcal{Y}^{1/2},
\end{eqnarray}
where $=_1$ comes from Eq.~\eqref{eq6-8-1:lma:thm 6.2}, and we use the following facts for $=_2$:
\begin{eqnarray}
\mathcal{Y}+\mathcal{A}_n &\rightarrow& \mathcal{Y},~~\mbox{in RTT as $\mathcal{A}_n \rightarrow \mathcal{O}$}\nonumber \\
(\mathcal{Y}+\mathcal{A}_n)^{-1/2}\mathcal{Y}^{1/2} &\rightarrow& \mbox{PR}(\mathcal{Y}),~~\mbox{in RTT as $\mathcal{A}_n \rightarrow \mathcal{O}$} \nonumber \\
\mathcal{Y}^{1/2}(\mathcal{Y}+\mathcal{A}_n)^{-1/2} &\rightarrow& \mbox{PR}(\mathcal{Y}).~~\mbox{in RTT as $\mathcal{A}_n \rightarrow \mathcal{O}$}
\end{eqnarray}
$\hfill \Box$

We define the following two sets for the pairs of tensors $(\mathcal{X}, \mathcal{Y})$ with respect to the order relation of tensors $\mathcal{X}$ and $\mathcal{Y}$. 
\begin{eqnarray}
\mathfrak{S}_{\preceq} &\define&\{(\mathcal{X},\mathcal{Y}): \mathcal{X} \preceq c \mathcal{Y} \mbox{~for some $c > 0$}\},\nonumber \\ 
\mathfrak{S}_{\succeq} &\define& \{(\mathcal{X},\mathcal{Y}): c\mathcal{X} \succeq \mathcal{Y} \mbox{~for some $c > 0$}\}. 
\end{eqnarray} 
Then, we have the following Theorem about the jointly tensor convexity on sets $\mathfrak{S}_{\preceq}$ and $\mathfrak{S}_{\succeq}$.
\begin{theorem}\label{thm:prop 6.5}
If $g \in \mbox{TC}$ and $g(0^+)=0$, then $(\mathcal{X},\mathcal{Y}) \rightarrow \mathcal{X} \#_g \mathcal{Y}$ is jointly tensor conex on $\mathfrak{S}_{\preceq}$. On the other hand, if $g'(\infty) < \infty$, then $(\mathcal{X},\mathcal{Y}) \rightarrow \mathcal{X} \#_g \mathcal{Y}$ is jointly tensor conex on $\mathfrak{S}_{\succeq}$. 
\end{theorem}
\textbf{Proof:}
For the case with $g(0^+) < \infty$,  $(\mathcal{X},\mathcal{Y}) \rightarrow \mathcal{X} \#_g \mathcal{Y}$ is jointly tensor conex on $\mathfrak{S}_{\preceq}$ from Lemma~\ref{lma:thm 6.2} and Theorem 2.2 in~\cite{ebadian2011perspectives}.

By defining $\tilde{g}(x) \define xg(x^{-1})$, we have the following relations $\tilde{g}'(\infty) = g(0^+)$ and $\mathcal{X} \#_{\tilde{g}} \mathcal{Y} = \mathcal{Y} \#_{g} \mathcal{X}$. By applying Lemma~\ref{lma:thm 6.2} and Theorem~\ref{thm:6.3} to $\tilde{g}(x)$, we have the following equivalent conditions:
\begin{enumerate}[label=(\roman*)]
\item For every PSD tensors $\mathcal{X}, \mathcal{Y}$ with $c\mathcal{X} \succeq \mathcal{Y}$ for some $c >0$, we have
\begin{eqnarray}
\lim\limits_{\epsilon \rightarrow 0^+} \left(\mathcal{Y}+\epsilon\mathcal{I}\right)\#_{\tilde{g}} \left(\mathcal{X}+\epsilon \mathcal{I}\right)&=&\lim\limits_{n \rightarrow \infty}\mathcal{Y}\#_{\tilde{g}} (\mathcal{X}+\mathcal{A}_n )\nonumber \\
&=& \mathcal{X}^{1/2}\star_N \tilde{g}(\eta(\mathcal{Y},\mathcal{X})) \star_N \mathcal{X}^{1/2},
\end{eqnarray}
where $\left\Vert \mathcal{A}_n \right\Vert_{\rho} \rightarrow \mathcal{O}$ for $n  \rightarrow \infty$;
\item $\tilde{g}'(\infty) < \infty$.
\end{enumerate}
For the case with $g'(\infty) < \infty$, $(\mathcal{X},\mathcal{Y}) \rightarrow \mathcal{X} \#_g \mathcal{Y}$ is jointly tensor conex on $\mathfrak{S}_{\succeq}$  from the aforementioned equivalent conditions and Theorem 2.2 in~\cite{ebadian2011perspectives} again since $\tilde{g}(x)$ is a convex function.
$\hfill \Box$

\subsection{Generalization Tail bounds for Bivariate Random Tensor Means from Part I}\label{sec:Generalization Tail bounds for Bivariate Random Tensor Means from Part I}

By assuming $\mathcal{Y}^{2^{k-1}} \preceq c \mathcal{X}^{2^{k-1}}$ for $k=1,2,\cdots,n$, we can represent $\mathcal{Z}_{k-1} \define \mathcal{X}^{-2^{k-2}} \mathcal{Y}^{2^{k-1}}\mathcal{X}^{-2^{k-2}}$ by $\eta(\mathcal{Y}^{2^{k-1}},\mathcal{X}^{2^{k-1}})$ according to the construction of the Eq.~\eqref{eq6-8-1:lma:thm 6.2} with respect to PSD tensors $\mathcal{X}$ and $\mathcal{Y}$. Therefore, we can extend those theorems and collaries presend in our Part I by allowing PSD random tensors based on the limitation method discussed in Section~\ref{sec:Limitation Method by Random Tensor Topology (RTT)}.  

Let us define the generalized product operation, denoted by $\acute{\prod}_{k=1}^n$, when the index upper bound is less than the index lower bound:
\begin{eqnarray}\label{eq4-1-2:thm:3.2}
\acute{\prod}_{k=1}^n a_i \define \begin{cases}
\prod_{k=1}^n a_i, ~\mbox{if $n \geq 1$};        \\
        1,~\mbox{if $n = 0$}.
        \end{cases}
\end{eqnarray}
where $a_i$ is the $i$-th real number.

Following Theorem is the extension of Theorem 4 in our Part I work.
\begin{theorem}\label{thm:3.2}
Given two random PSD tensors $\mathcal{X} $$\in$$ \mathbb{C}^{I_1 \times \cdots \times I_N \times I_1 \times \cdots \times I_N}$, $\mathcal{Y} $$\in$$ \mathbb{C}^{I_1 \times \cdots \times I_N \times I_1 \times \cdots \times I_N}$ and a PD determinstic tensor $\mathcal{C}$, if $q=2^n q_0 \geq 1$ with $1 \leq q_0 \leq 2$ and $n \in \mathbb{N}$, we set $\mathcal{Z}_{k-1} \define \eta(\mathcal{Y}^{2^{k-1}},\mathcal{X}^{2^{k-1}})$ by assuming that $\mathcal{Y}^{2^{k-1}} \preceq c \mathcal{X}^{2^{k-1}}$ for $k=1,2,\cdots,n$. We also assume that $\mathcal{X} $$\#_f$$ \mathcal{Y} $$\succeq$$ \mathcal{I}$ almost surely with $f $$\in$$\mbox{TMI}^{1}$, we have 
\begin{eqnarray}\label{eq1-1:thm:3.2}
\mathrm{Pr}\left(\mathcal{X}^q \#_f \mathcal{Y}^q \npreceq \mathcal{C} \right)
&\leq& \mathrm{Tr}\left(\mathbb{E}\left[\left(\Psi_{upper}\left(q,f,\mathcal{X},\mathcal{Y}\right)\lambda^{q-1}_{\min}\left(\mathcal{X}\#_f\mathcal{Y}\right)\mathcal{X}\#_f\mathcal{Y}\right)^p\right]\star_N \mathcal{C}^{-1} \right),
\end{eqnarray}
and
\begin{eqnarray}\label{eq1-2:thm:3.2}
\mathrm{Pr}\left(\Psi_{lower}\left(q,f,\mathcal{X},\mathcal{Y}\right)\lambda^{q-1}_{\max}\left(\mathcal{X}\#_f\mathcal{Y}\right)\mathcal{X}\#_f\mathcal{Y} \npreceq \mathcal{C} \right)
&\leq& \mathrm{Tr}\left(\mathbb{E}\left[\left(\mathcal{X}^q \#_f \mathcal{Y}^q\right)^p\right]\star_N \mathcal{C}^{-1} \right),
\end{eqnarray}
where $\Psi_{lower}\left(q,f,\mathcal{X},\mathcal{Y}\right)$ and $\Psi_{upper}\left(q,f,\mathcal{X},\mathcal{Y}\right)$ are two positive numbers defined by
\begin{eqnarray}
\Psi_{lower}\left(q,f,\mathcal{X},\mathcal{Y}\right)&\define&\lambda_{\min}\left(f^{-q_0}\left(\mathcal{Z}_n\right)f\left(\mathcal{Z}^{q_0}_n\right)\right) \acute{\prod}_{k=1}^n \lambda_{\min}\left(f^{-2}\left(\mathcal{Z}_{k-1}\right)f\left(\mathcal{Z}_{k-1}^{2}\right)\right) \nonumber \\
\Psi_{upper}\left(q,f,\mathcal{X},\mathcal{Y}\right)&\define&\lambda_{\max}\left(f^{-q_0}\left(\mathcal{Z}_n\right)(f\left(\mathcal{Z}_n^{q_0}\right)\right)\acute{\prod}_{k=1}^n \lambda_{\max}\left(f^{-2}\left(\mathcal{Z}_{k-1}\right)f\left(\mathcal{Z}_{k-1}^{2}\right)\right).
\end{eqnarray}
Note that the definition of $\acute{\prod}$ is provided by Eq.~\eqref{eq4-1-2:thm:3.2}.

For $0 < q  \leq1$, we have 
\begin{eqnarray}\label{eq2-1:thm:3.2}
\mathrm{Pr}\left(\mathcal{X}^q \#_f \mathcal{Y}^q \npreceq \mathcal{C} \right)
&\leq& \mathrm{Tr}\left(\mathbb{E}\left[\left(\lambda_{\min}\left(f^{-q}\left(\mathcal{Z}_0\right)(f\left(\mathcal{Z}_0^q\right)\right)\lambda^{q-1}_{\min}\left(\mathcal{X}\#_f\mathcal{Y}\right)\mathcal{X}\#_f\mathcal{Y}\right)^p\right]\star_N \mathcal{C}^{-1} \right),
\end{eqnarray}
and
\begin{eqnarray}\label{eq2-2:thm:3.2}
\mathrm{Pr}\left( \lambda_{\max}\left(f^{-q}\left(\mathcal{Z}_0\right)f\left(\mathcal{Z}_0^q\right)\right)\lambda^{q-1}_{\max}\left(\mathcal{X}\#_f\mathcal{Y}\right)\mathcal{X}\#_f\mathcal{Y} \npreceq \mathcal{C} \right)
&\leq& \mathrm{Tr}\left(\mathbb{E}\left[\left(\mathcal{X}^q \#_f \mathcal{Y}^q\right)^p\right]\star_N \mathcal{C}^{-1} \right).
\end{eqnarray}
where $p \geq 1$. 
\end{theorem}

Following Theorem is the extension of Theorem 5 in our Part I work.
\begin{theorem}\label{thm:Prop. 3.9}
Given two random PSD tensors $\mathcal{X},\mathcal{Y} $$\in$$ \mathbb{C}^{I_1 \times \cdots \times I_N \times I_1 \times \cdots \times I_N}$, and a PD determinstic tensor $\mathcal{C}$, if $q=2^n q_0 \geq 1$ with $1 \leq q_0 \leq 2$, we set $\mathcal{Z}_{k-1} \define \eta(\mathcal{Y}^{2^{k-1}},\mathcal{X}^{2^{k-1}})$ by assuming that $\mathcal{Y}^{2^{k-1}} \preceq c \mathcal{X}^{2^{k-1}}$ for $k=1,2,\cdots,n$. We assume that $\mathcal{X} $$\#_h$$ \mathcal{Y} $$\preceq$$ \mathcal{I}$ almost surely with $h $$\in$$\mbox{TMD}^{1}$, we have 
\begin{eqnarray}\label{eq1-1:Prop. 3.9}
\mathrm{Pr}\left(\mathcal{X}^q \#_h \mathcal{Y}^q \npreceq \mathcal{C} \right)
&\leq& \mathrm{Tr}\left(\mathbb{E}\left[\left(\Phi_{upper}\left(q,h,\mathcal{X},\mathcal{Y}\right)\lambda^{q-1}_{\min}\left(\mathcal{X}\#_h\mathcal{Y}\right)\mathcal{X}\#_h\mathcal{Y}\right)^p\right]\star_N \mathcal{C}^{-1} \right),
\end{eqnarray}
and
\begin{eqnarray}\label{eq1-2:Prop. 3.9}
\mathrm{Pr}\left(\Phi_{lower}\left(q,h,\mathcal{X},\mathcal{Y}\right)\lambda^{q-1}_{\max}\left(\mathcal{X}\#_h\mathcal{Y}\right)\mathcal{X}\#_h\mathcal{Y} \npreceq \mathcal{C} \right)
&\leq& \mathrm{Tr}\left(\mathbb{E}\left[\left(\mathcal{X}^q \#_h \mathcal{Y}^q\right)^p\right]\star_N \mathcal{C}^{-1} \right),
\end{eqnarray}
where $\Phi_{lower}\left(q,h,\mathcal{X},\mathcal{Y}\right)$ and $\Phi_{upper}\left(q,h,\mathcal{X},\mathcal{Y}\right)$ are two positive numbers defined by
\begin{eqnarray}
\Phi_{lower}\left(q,h,\mathcal{X},\mathcal{Y}\right)&\define&\lambda_{\min}\left(h^{-q_0}\left(\mathcal{Z}_n\right)h\left(\mathcal{Z}^{q_0}_n\right)\right) \acute{\prod}_{k=1}^n \lambda_{\min}\left(h^{-2}\left(\mathcal{Z}_{k-1}\right)h\left(\mathcal{Z}_{k-1}^{2}\right)\right) \nonumber \\
\Phi_{upper}\left(q,h,\mathcal{X},\mathcal{Y}\right)&\define&\lambda_{\max}\left(h^{-q_0}\left(\mathcal{Z}_n\right)(h\left(\mathcal{Z}_n^{q_0}\right)\right)\acute{\prod}_{k=1}^n \lambda_{\max}\left(h^{-2}\left(\mathcal{Z}_{k-1}\right)h\left(\mathcal{Z}_{k-1}^{2}\right)\right).
\end{eqnarray}
Note that the definition of $\acute{\prod}$ is provided by Eq.~\eqref{eq4-1-2:thm:3.2}.

For $0 < q  \leq1$, we have 
\begin{eqnarray}\label{eq2-1:Prop. 3.9}
\mathrm{Pr}\left(\mathcal{X}^q \#_h \mathcal{Y}^q \npreceq \mathcal{C} \right)
&\leq& \mathrm{Tr}\left(\mathbb{E}\left[\left(\lambda_{\max}\left(h^{-q}\left(\mathcal{Z}_0\right)(h\left(\mathcal{Z}_0^q\right)\right)\lambda^{q-1}_{\min}\left(\mathcal{X}\#_h\mathcal{Y}\right)\mathcal{X}\#_h\mathcal{Y}\right)^p\right]\star_N \mathcal{C}^{-1} \right),
\end{eqnarray}
and
\begin{eqnarray}\label{eq2-2:Prop. 3.9}
\mathrm{Pr}\left( \lambda_{\min}\left(h^{-q}\left(\mathcal{Z}_0\right)h\left(\mathcal{Z}_0^q\right)\right)\lambda^{q-1}_{\max}\left(\mathcal{X}\#_h\mathcal{Y}\right)\mathcal{X}\#_h\mathcal{Y} \npreceq \mathcal{C} \right)
&\leq& \mathrm{Tr}\left(\mathbb{E}\left[\left(\mathcal{X}^q \#_h \mathcal{Y}^q\right)^p\right]\star_N \mathcal{C}^{-1} \right).
\end{eqnarray}
where $p \geq 1$. 
\end{theorem}

Following Theorem is the extension of Theorem 6 in our Part I work.
\begin{theorem}\label{thm:Prop. 3.10}
Given a PD tensor $\mathcal{X} $$\in$$ \mathbb{C}^{I_1 \times \cdots \times I_N \times I_1 \times \cdots \times I_N}$, a PSD tensor $\mathcal{Y} $$\in$$ \mathbb{C}^{I_1 \times \cdots \times I_N \times I_1 \times \cdots \times I_N}$, and a PD determinstic tensor $\mathcal{C} $$\in$$ \mathbb{C}^{I_1 \times \cdots \times I_N \times I_1 \times \cdots \times I_N}$, we will set $\mathcal{Z} \define \eta^{-1}(\mathcal{Y},\mathcal{X})$ by assuming that $\mathcal{Y} \preceq c \mathcal{X}$ and the tensor $\eta(\mathcal{Y},\mathcal{X})$ being invertible. Let $g $$\in$$\mbox{TC}^{1}$, if $\mathcal{X} \#_g \mathcal{Y} \preceq \mathcal{I}$ almost surely, and $p,q \geq 1$, we have 
\begin{eqnarray}\label{eq1-1:Prop. 3.10}
\mathrm{Pr}\left(\mathcal{X}^q \#_g \mathcal{Y}^q \npreceq \mathcal{C}\right) \leq
\mathrm{Tr}\left(\mathbb{E}\left[\left(\mathrm{K}_1 \lambda^{1-q}_{\min}\left(\mathcal{X}\#_{g}\mathcal{Y}\right)
\lambda_{\max}\left(g^{-q}(\mathcal{Z})g(\mathcal{Z}^q)\right)\mathrm{K}_2\mathcal{I}\right)^p\right]\star_N \mathcal{C}^{-1}\right)
\end{eqnarray}
where $\mathrm{K}_1$ and $\mathrm{K}_2$ are set as
\begin{eqnarray}
\mathrm{K}_1 &\define& \mathrm{K}\left(\lambda^{-1}_{\max}\left(\mathcal{X}\right),\lambda^{-1}_{\min}\left(\mathcal{X}\right),q-1\right) \nonumber \\
\mathrm{K}_2 &\define& \mathrm{K}\left(\lambda^{-1}_{\max}\left(\mathcal{X}\right),\lambda^{-1}_{\min}\left(\mathcal{X}\right),2q-1\right).
\end{eqnarray}
Moreover, if $\mathcal{X} \#_g \mathcal{Y} \succeq \mathcal{I}$ almost surely, we have 
\begin{eqnarray}\label{eq1-2:Prop. 3.10}
\mathrm{Pr}\left(\lambda^{1-q}_{\min}\left(\mathcal{X}\#_{g}\mathcal{Y}\right)
\lambda_{\max}\left(g^{-q}(\mathcal{Z})g(\mathcal{Z}^q)\right) \mathrm{K}_2^{-1}\mathcal{I}\npreceq \mathcal{C}\right) \leq \mathrm{Tr}\left(\mathbb{E}\left[\left(\mathcal{X}^q \#_g \mathcal{Y}^q\right)^p\right] \star_N \mathcal{C}^{-1}\right)
\end{eqnarray}
\end{theorem}

Following Corollary is the extension of Corollary 3 in our Part I work.
\begin{corollary}\label{thm:3.2 major}
Given two random PSD tensors $\mathcal{X} $$\in$$ \mathbb{C}^{I_1 \times \cdots \times I_N \times I_1 \times \cdots \times I_N}$, $\mathcal{Y} $$\in$$ \mathbb{C}^{I_1 \times \cdots \times I_N \times I_1 \times \cdots \times I_N}$ and a PD determinstic tensor $\mathcal{C}$, if $q=2^n q_0 \geq 1$ with $1 \leq q_0 \leq 2$, we set $\mathcal{Z}_{k-1} =\eta(\mathcal{Y}^{2^{k-1}},\mathcal{X}^{2^{k-1}})$ by assuming that $\mathcal{Y}^{2^{k-1}} \preceq c \mathcal{X}^{2^{k-1}}$ for $k=1,2,\cdots,n$. We assume that $\mathcal{X} $$\#_f$$ \mathcal{Y} $$\succeq$$ \mathcal{I}$ almost surely with $f $$\in$$\mbox{TMI}^{1}$. 
Then, we have 
\begin{eqnarray}\label{eq1-1:thm:3.2 major}
\mathrm{Pr}\left(\sum\limits_{i=1}^{k}\lambda_i\left(\Psi_{lower}\left(q,f,\mathcal{X},\mathcal{Y}\right)\lambda^{q-1}_{\max}\left(\mathcal{X}\#_f\mathcal{Y}\right)\mathcal{X}\#_f\mathcal{Y}\right) \geq \kappa \right)
\leq
\mathrm{Pr}\left(\sum\limits_{i=1}^{k}\lambda_i\left(\mathcal{X}^q \#_f \mathcal{Y}^q \right) \geq \kappa \right) \nonumber \\
\leq \mathrm{Pr}\left(\sum\limits_{i=1}^{k}\lambda_i\left(\Psi_{upper}\left(q,f,\mathcal{X},\mathcal{Y}\right)\lambda^{q-1}_{\min}\left(\mathcal{X}\#_f\mathcal{Y}\right)\mathcal{X}\#_f\mathcal{Y}\right) \geq \kappa \right),~~~~~~~~~~~~~~~~~~~~~~~~~~~~~~~~~~~~~~~
\end{eqnarray}
and
\begin{eqnarray}\label{eq1-2:thm:3.2 major}
\mathrm{Pr}\left(\prod\limits_{i=1}^{k}\lambda_i\left(\Psi_{lower}\left(q,f,\mathcal{X},\mathcal{Y}\right)\lambda^{q-1}_{\max}\left(\mathcal{X}\#_f\mathcal{Y}\right)\mathcal{X}\#_f\mathcal{Y}\right) \geq \kappa \right)
\leq
\mathrm{Pr}\left(\prod\limits_{i=1}^{k}\lambda_i\left(\mathcal{X}^q \#_f \mathcal{Y}^q \right) \geq \kappa \right) \nonumber \\
\leq \mathrm{Pr}\left(\prod\limits_{i=1}^{k}\lambda_i\left(\Psi_{upper}\left(q,f,\mathcal{X},\mathcal{Y}\right)\lambda^{q-1}_{\min}\left(\mathcal{X}\#_f\mathcal{Y}\right)\mathcal{X}\#_f\mathcal{Y}\right) \geq \kappa \right).~~~~~~~~~~~~~~~~~~~~~~~~~~~~~~~~~~~~~~~
\end{eqnarray}

For $0 < q  \leq1$, we have 
\begin{eqnarray}\label{eq2-1:thm:3.2 major}
\mathrm{Pr}\left(\sum\limits_{i=1}^k\lambda_i\left(\lambda_{\max}\left(f^{-q}\left(\mathcal{Z}_0\right)f\left(\mathcal{Z}_0^q\right)\right)\lambda^{q-1}_{\max}\left(\mathcal{X}\#_f\mathcal{Y}\right)\mathcal{X}\#_f\mathcal{Y}\right) \geq \kappa \right) \leq \mathrm{Pr}\left(\sum\limits_{i=1}^k\lambda_i\left(\mathcal{X}^q \#_f \mathcal{Y}^q \right) \geq \kappa \right) \nonumber \\
\leq  \mathrm{Pr}\left(\sum\limits_{i=1}^k\lambda_i\left(\lambda_{\min}\left(f^{-q}\left(\mathcal{Z}_0\right)(f\left(\mathcal{Z}_0^q\right)\right)\lambda^{q-1}_{\min}\left(\mathcal{X}\#_f\mathcal{Y}\right)\mathcal{X}\#_f\mathcal{Y} \right) \geq \kappa \right),~~~~~~~~~~~~~~~~~~~~~~~~~~~~~~~~~~~~~
\end{eqnarray}
and
\begin{eqnarray}\label{eq2-2:thm:3.2 major}
\mathrm{Pr}\left(\prod\limits_{i=1}^k\lambda_i\left(\lambda_{\max}\left(f^{-q}\left(\mathcal{Z}_0\right)f\left(\mathcal{Z}_0^q\right)\right)\lambda^{q-1}_{\max}\left(\mathcal{X}\#_f\mathcal{Y}\right)\mathcal{X}\#_f\mathcal{Y}\right) \geq \kappa \right) \leq \mathrm{Pr}\left(\prod\limits_{i=1}^k\lambda_i\left(\mathcal{X}^q \#_f \mathcal{Y}^q \right) \geq \kappa \right) \nonumber \\
\leq  \mathrm{Pr}\left(\prod\limits_{i=1}^k\lambda_i\left(\lambda_{\min}\left(f^{-q}\left(\mathcal{Z}_0\right)(f\left(\mathcal{Z}_0^q\right)\right)\lambda^{q-1}_{\min}\left(\mathcal{X}\#_f\mathcal{Y}\right)\mathcal{X}\#_f\mathcal{Y} \right) \geq \kappa \right).~~~~~~~~~~~~~~~~~~~~~~~~~~~~~~~~~~~~~
\end{eqnarray}
\end{corollary}

Following Corollary is the extension of Corollary 4 in our Part I work.
\begin{corollary}\label{Prop. 3.9 major}
Given two random PSD tensors $\mathcal{X}, \mathcal{Y} $$\in$$ \mathbb{C}^{I_1 \times \cdots \times I_N \times I_1 \times \cdots \times I_N}$, and a PD determinstic tensor $\mathcal{C}$, if $q=2^n q_0 \geq 1$ with $1 \leq q_0 \leq 2$, we set $\mathcal{Z}_{k-1}=\eta(\mathcal{Y}^{2^{k-1}},\mathcal{X}^{2^{k-1}})$ by assuming that $\mathcal{Y}^{2^{k-1}} \preceq c \mathcal{X}^{2^{k-1}}$ for $k=1,2,\cdots,n$. We assume that $\mathcal{X} $$\#_h$$ \mathcal{Y} $$\preceq$$ \mathcal{I}$ almost surely with $h $$\in$$\mbox{TMD}^{1}$. 
Then, we have 
\begin{eqnarray}\label{eq1-1:Prop. 3.9 major}
\mathrm{Pr}\left(\sum\limits_{i=1}^k\lambda_i\left(\Phi_{lower}\left(q,h,\mathcal{X},\mathcal{Y}\right)\lambda^{q-1}_{\max}\left(\mathcal{X}\#_h\mathcal{Y}\right)\mathcal{X}\#_h\mathcal{Y}\right) \geq \kappa \right) \leq \mathrm{Pr}\left(\sum\limits_{i=1}^k\lambda_i\left(\mathcal{X}^q \#_h \mathcal{Y}^q \right) \geq \kappa \right)  \nonumber \\
\leq \mathrm{Pr}\left(\sum\limits_{i=1}^k\lambda_i\left(\Phi_{upper}\left(q,h,\mathcal{X},\mathcal{Y}\right)\lambda^{q-1}_{\min}\left(\mathcal{X}\#_h\mathcal{Y}\right)\mathcal{X}\#_h\mathcal{Y}\right) \geq \kappa \right),~~~~~~~~~~~~~~~~~~~~~~~~~~~~~~~~~~~~~~
\end{eqnarray}
and
\begin{eqnarray}\label{eq1-2:Prop. 3.9 major}
\mathrm{Pr}\left(\prod\limits_{i=1}^k\lambda_i\left(\Phi_{lower}\left(q,h,\mathcal{X},\mathcal{Y}\right)\lambda^{q-1}_{\max}\left(\mathcal{X}\#_h\mathcal{Y}\right)\mathcal{X}\#_h\mathcal{Y}\right) \geq \kappa \right) \leq \mathrm{Pr}\left(\prod\limits_{i=1}^k\lambda_i\left(\mathcal{X}^q \#_h \mathcal{Y}^q \right) \geq \kappa \right)  \nonumber \\
\leq \mathrm{Pr}\left(\prod\limits_{i=1}^k\lambda_i\left(\Phi_{upper}\left(q,h,\mathcal{X},\mathcal{Y}\right)\lambda^{q-1}_{\min}\left(\mathcal{X}\#_h\mathcal{Y}\right)\mathcal{X}\#_h\mathcal{Y}\right) \geq \kappa \right).~~~~~~~~~~~~~~~~~~~~~~~~~~~~~~~~~~~~~
\end{eqnarray}

For $0 < q  \leq1$, we have 
\begin{eqnarray}\label{eq2-1:Prop. 3.9 major}
\mathrm{Pr}\left(\sum\limits_{i=1}^k\lambda_i\left(\lambda_{\min}\left(h^{-q}\left(\mathcal{Z}_0\right)h\left(\mathcal{Z}_0^q\right)\right)\lambda^{q-1}_{\max}\left(\mathcal{X}\#_h\mathcal{Y}\right)\mathcal{X}\#_h\mathcal{Y}\right) \geq \kappa \right) \leq \mathrm{Pr}\left(\sum\limits_{i=1}^k\lambda_i\left(\mathcal{X}^q \#_h \mathcal{Y}^q \right) \geq \kappa \right) \nonumber \\
\leq  \mathrm{Pr}\left(\sum\limits_{i=1}^k\lambda_i\left(\lambda_{\max}\left(h^{-q}\left(\mathcal{Z}_0\right)(h\left(\mathcal{Z}_0^q\right)\right)\lambda^{q-1}_{\min}\left(\mathcal{X}\#_h\mathcal{Y}\right)\mathcal{X}\#_h\mathcal{Y}\right) \geq \kappa \right),~~~~~~~~~~~~~~~~~~~~~~~~~~~~~~~~~~~~~
\end{eqnarray}
and
\begin{eqnarray}\label{eq2-2:Prop. 3.9 major}
\mathrm{Pr}\left(\prod\limits_{i=1}^k\lambda_i\left(\lambda_{\min}\left(h^{-q}\left(\mathcal{Z}_0\right)h\left(\mathcal{Z}_0^q\right)\right)\lambda^{q-1}_{\max}\left(\mathcal{X}\#_h\mathcal{Y}\right)\mathcal{X}\#_h\mathcal{Y}\right) \geq \kappa \right) \leq \mathrm{Pr}\left(\prod\limits_{i=1}^k\lambda_i\left(\mathcal{X}^q \#_h \mathcal{Y}^q \right) \geq \kappa \right) \nonumber \\
\leq  \mathrm{Pr}\left(\prod\limits_{i=1}^k\lambda_i\left(\lambda_{\max}\left(h^{-q}\left(\mathcal{Z}_0\right)(h\left(\mathcal{Z}_0^q\right)\right)\lambda^{q-1}_{\min}\left(\mathcal{X}\#_h\mathcal{Y}\right)\mathcal{X}\#_h\mathcal{Y}\right) \geq \kappa \right).~~~~~~~~~~~~~~~~~~~~~~~~~~~~~~~~~~~~~
\end{eqnarray}
\end{corollary}

Following Corollary is the extension of Corollary 5 in our Part I work.
\begin{corollary}\label{thm:Prop. 3.10 major}
Given a PD random tensors $\mathcal{X} $$\in$$ \mathbb{C}^{I_1 \times \cdots \times I_N \times I_1 \times \cdots \times I_N}$ and a PSD random tensors $\mathcal{Y} $$\in$$ \mathbb{C}^{I_1 \times \cdots \times I_N \times I_1 \times \cdots \times I_N}$, and a PD determinstic tensor $\mathcal{C} $$\in$$ \mathbb{C}^{I_1 \times \cdots \times I_N \times I_1 \times \cdots \times I_N}$, we will set $\mathcal{Z}=\eta^{-1}(\mathcal{Y},\mathcal{X})$ by assuming that $\mathcal{Y}\preceq c \mathcal{X}$ and the invertibility of the tensor $\eta(\mathcal{Y},\mathcal{X})$. Let $g $$\in$$\mbox{TC}^{1}$, if $\mathcal{X} \#_g \mathcal{Y} \preceq \mathcal{I}$ almost surely and $q \geq 1$, we have 
\begin{eqnarray}\label{eq1-1:thm:Prop. 3.10 major}
\mathrm{Pr}\left(\sum\limits_{i=1}^k\lambda_i\left(\mathcal{X}^q \#_g \mathcal{Y}^q  \right) \geq \kappa \right) &\leq& \mathrm{Pr}\Bigg(\sum\limits_{i=1}^k\lambda_i\Big(\mathrm{K}\left(\lambda^{-1}_{\max}\left(\mathcal{X}\right),\lambda^{-1}_{\min}\left(\mathcal{X}\right),q-1\right)\lambda^{1-q}_{\min}\left(\mathcal{X}\#_{g}\mathcal{Y}\right) \nonumber \\
&  & \lambda_{\max}\left(g^{-q}(\mathcal{Z})g(\mathcal{Z}^q)\right)\mathrm{K}\left(\lambda^{-1}_{\max}\left(\mathcal{X}\right),\lambda^{-1}_{\min}\left(\mathcal{X}\right),2q-1\right)\mathcal{I}\Big) \geq \kappa \Bigg),\nonumber \\
\end{eqnarray}
and
\begin{eqnarray}\label{eq1-2:thm:Prop. 3.10 major}
\mathrm{Pr}\left(\prod\limits_{i=1}^k\lambda_i\left(\mathcal{X}^q \#_g \mathcal{Y}^q  \right) \geq \kappa \right) &\leq& \mathrm{Pr}\Bigg(\prod\limits_{i=1}^k\lambda_i\Big(\mathrm{K}\left(\lambda^{-1}_{\max}\left(\mathcal{X}\right),\lambda^{-1}_{\min}\left(\mathcal{X}\right),q-1\right)\lambda^{1-q}_{\min}\left(\mathcal{X}\#_{g}\mathcal{Y}\right) \nonumber \\
&  & \lambda_{\max}\left(g^{-q}(\mathcal{Z})g(\mathcal{Z}^q)\right)\mathrm{K}\left(\lambda^{-1}_{\max}\left(\mathcal{X}\right),\lambda^{-1}_{\min}\left(\mathcal{X}\right),2q-1\right)\mathcal{I}\Big) \geq \kappa \Bigg)\nonumber \\
\end{eqnarray}

Moreover, if $\mathcal{X} \#_g \mathcal{Y} \succeq \mathcal{I}$ almost surely, we have 
\begin{eqnarray}\label{eq2-1:thm:Prop. 3.10 major}
\mathrm{Pr}\left(\sum\limits_{i=1}^k\lambda_i\left(\mathcal{X}^q \#_g \mathcal{Y}^q  \right) \geq \kappa \right) &\geq& \mathrm{Pr}\Bigg(\sum\limits_{i=1}^k\lambda_i\Big(\lambda^{1-q}_{\min}\left(\mathcal{X}\#_{g}\mathcal{Y}\right)\lambda_{\max}\left(g^{-q}(\mathcal{Z})g(\mathcal{Z}^q)\right) \nonumber \\
&  & \mathrm{K}^{-1}\left(\lambda^{-1}_{\max}\left(\mathcal{X}\right),\lambda^{-1}_{\min}\left(\mathcal{X}\right),2q-1\right)\mathcal{I}\Big) \geq \kappa \Bigg),
\end{eqnarray}
and
\begin{eqnarray}\label{eq2-2:thm:Prop. 3.10 major}
\mathrm{Pr}\left(\prod\limits_{i=1}^k\lambda_i\left(\mathcal{X}^q \#_g \mathcal{Y}^q  \right) \geq \kappa \right) &\geq& \mathrm{Pr}\Bigg(\prod\limits_{i=1}^k\lambda_i\Big(\lambda^{1-q}_{\min}\left(\mathcal{X}\#_{g}\mathcal{Y}\right)\lambda_{\max}\left(g^{-q}(\mathcal{Z})g(\mathcal{Z}^q)\right) \nonumber \\
&  & \mathrm{K}^{-1}\left(\lambda^{-1}_{\max}\left(\mathcal{X}\right),\lambda^{-1}_{\min}\left(\mathcal{X}\right),2q-1\right)\mathcal{I}\Big) \geq \kappa \Bigg).
\end{eqnarray}
\end{corollary}

\section{Applications: Tensor Data Processing}\label{sec:Applications: Tensor Data Processing} 

In this section, we will derive two new tail bounds for bivariate tensor means under two tensor data processing methods: data fusion and linear transformation.

In Theorem~\ref{thm:tensor processing inequality}, we will study the relationship for tensor data fusion via addition before taking the tensor mean operation, i.e., $(\mathcal{X}_1+\mathcal{X}_2) \#_g (\mathcal{Y}_1+\mathcal{Y}_2)$, and after taking the tensor mean operation, i.e., $\mathcal{X}_1 \#_g \mathcal{Y}_1 + \mathcal{X}_2 \#_g \mathcal{Y}_2$.

\begin{theorem}[Tensor Mean Fusion Inequality]\label{thm:tensor processing inequality}
Given $g$ be an operator convex function on $(0, \infty)$ with $g(0^+) < \infty$, and four random PSD tensors $\mathcal{X}_1, \mathcal{Y}_1, \mathcal{X}_2, \mathcal{Y}_2 \in \mathbb{C}^{I_1 \times \dots \times I_N \times I_1 \times \dots \times I_N}$ with $(\mathcal{X}_1, \mathcal{Y}_1), (\mathcal{X}_2, \mathcal{Y}_2) \in \mathfrak{S}_{\preceq}$, we have 
\begin{eqnarray}\label{eq1-1:thm:tensor processing inequality}
(\mathcal{X}_1+\mathcal{X}_2) \#_g (\mathcal{Y}_1+\mathcal{Y}_2)
\preceq \mathcal{X}_1 \#_g \mathcal{Y}_1 + \mathcal{X}_2 \#_g \mathcal{Y}_2,
\end{eqnarray}
and, given a PD tensor $\mathcal{C} \in \mathbb{C}^{I_1 \times \dots \times I_N \times I_1 \times \dots \times I_N}$, we have
\begin{eqnarray}\label{eq1-2:thm:tensor processing inequality}
\mathrm{Pr}\left((\mathcal{X}_1+\mathcal{X}_2) \#_g (\mathcal{Y}_1+\mathcal{Y}_2)\npreceq \mathcal{C}\right)
\leq \mathrm{Tr}\left(\mathbb{E}\left[\left(\mathcal{X}_1 \#_g \mathcal{Y}_1 + \mathcal{X}_2 \#_g \mathcal{Y}_2\right)^q\right] \star_N \mathcal{C}^{-1}\right),
\end{eqnarray}
where $q \geq 1$.
\end{theorem}
\textbf{Proof:}

From Lemma~\ref{lma:prop 2.2} and Theorem~\ref{thm:prop 6.5}, we have 
\begin{eqnarray}
(\lambda \mathcal{X}_1+\overline{\lambda}\mathcal{X}_2) \#_g (\lambda \mathcal{Y}_1+\overline{\lambda}\mathcal{Y}_2)
\preceq \lambda \mathcal{X}_1 \#_g \mathcal{Y}_1 + \overline{\lambda}\mathcal{X}_2 \#_g \mathcal{Y}_2,
\end{eqnarray}
then, we have Eq.~\eqref{eq1-1:thm:tensor processing inequality} by setting $\lambda=1/2$ and using the property $(\lambda\mathcal{X}) \#_g (\lambda\mathcal{Y}) = \lambda(\mathcal{X} \#_g \mathcal{Y})$ for any tensors $\mathcal{X}, \mathcal{Y}$. Eq.~\eqref{eq1-2:thm:tensor processing inequality} is obtained by applying Lemma~\ref{lma:Loewner ordering with Markov Cheb inequalities} to Eq.~\eqref{eq1-1:thm:tensor processing inequality}.

Note that Eq.~\eqref{eq1-1:thm:tensor processing inequality} and Eq.~\eqref{eq1-2:thm:tensor processing inequality} also valid for $(\mathcal{X}_1, \mathcal{Y}_1), (\mathcal{X}_2, \mathcal{Y}_2) \in \mathfrak{S}_{\succeq}$.
$\hfill \Box$

Another important tail bound for bivariate tensor mean inequality is the monotonicity under positive linear transform for tensor data.
\begin{theorem}[Tensor Mean Inequality Under Linear Transform]\label{thm:Tensor Mean Inequality Under Linear Transform}
Given $g$ be an operator convex function on $(0, \infty)$ with $g(0^+) < \infty$, two random PSD tensors $\mathcal{X}, \mathcal{Y} \in \mathbb{C}^{I_1 \times \dots \times I_N \times I_1 \times \dots \times I_N}$ with $(\mathcal{X}, \mathcal{Y}) \in \mathfrak{S}_{\preceq}$ and linear transformation $\mathfrak{L}$ between tensors, we have 
\begin{eqnarray}\label{eq1-1:thm:Tensor Mean Inequality Under Linear Transform}
\mathfrak{L}\left(\mathcal{X}\#_{g}\mathcal{Y}\right)&\succeq& \mathfrak{L}(\mathcal{X}) \#_{g} \mathfrak{L}(\mathcal{Y}),
\end{eqnarray}
and, given a PD tensor $\mathcal{C} \in \mathbb{C}^{I_1 \times \dots \times I_N \times I_1 \times \dots \times I_N}$, we have
\begin{eqnarray}\label{eq1-2:thm:Tensor Mean Inequality Under Linear Transform}
\mathrm{Pr}\left(\mathfrak{L}(\mathcal{X}) \#_{g} \mathfrak{L}(\mathcal{Y}) \npreceq \mathcal{C}\right)
\leq \mathrm{Tr}\left(\mathbb{E}\left[\mathfrak{L}^q\left(\mathcal{X} \#_{g} \mathcal{Y}\right)\right] \star_N \mathcal{C}^{-1}\right),
\end{eqnarray}
where $q \geq 1$.
\end{theorem}
\textbf{Proof:}
We begin by assuming that both random tensors $\mathcal{X}$ and $\mathcal{Y}$ are PD tensors. Then, we have
\begin{eqnarray}\label{eq2:thm:Tensor Mean Inequality Under Linear Transform}
\mathfrak{L}(\mathcal{X}\#_{g}\mathcal{Y})&=&\mathfrak{L}^{1/2}(\mathcal{Y})\left[\mathfrak{L}^{-1/2}(\mathcal{Y})\mathfrak{L}\left(\mathcal{Y}^{1/2}g(\mathcal{Y}^{-1/2}\mathcal{X}\mathcal{Y}^{-1/2})\mathcal{Y}^{1/2}\right)\mathfrak{L}^{-1/2}(\mathcal{Y})\right]\mathfrak{L}^{1/2}(\mathcal{Y}) \nonumber \\
&\succeq_1& \mathfrak{L}^{1/2}(\mathcal{Y})g\left(\mathfrak{L}^{-1/2}(\mathcal{Y})\mathfrak{L}\left(\mathcal{X}\right)\mathfrak{L}^{-1/2}(\mathcal{Y})\right)\mathfrak{L}^{1/2}(\mathcal{Y}) \nonumber \\
&=&  \mathfrak{L}(\mathcal{X}) \#_{g} \mathfrak{L}(\mathcal{Y}),
\end{eqnarray}
where $\succeq_1$ comes from operator Jensen inequality.

If both random tensors $\mathcal{X}$ and $\mathcal{Y}$ are PSD tensors, we can approximate $\mathcal{X}$ and $\mathcal{Y}$ by PD tensors as $(\mathcal{X}+\epsilon \mathcal{I})$ and $(\mathcal{Y}+\epsilon \mathcal{I})$. Therefore, we have
\begin{eqnarray}\label{eq4:thm:Tensor Mean Inequality Under Linear Transform}
\mathfrak{L}\left((\mathcal{X}+\epsilon \mathcal{I})\#_{g}(\mathcal{Y}+\epsilon \mathcal{I})\right)&\succeq& \mathfrak{L}(\mathcal{X}+\epsilon \mathcal{I})\#_{g}\mathfrak{L}(\mathcal{Y}+\epsilon \mathcal{I}),
\end{eqnarray}
and, Eq.~\eqref{eq1-1:thm:Tensor Mean Inequality Under Linear Transform} is established from Lemma~\ref{lma:thm 6.2}. 

Eq.~\eqref{eq1-2:thm:Tensor Mean Inequality Under Linear Transform} is obtained by applying Lemma~\ref{lma:Loewner ordering with Markov Cheb inequalities} to Eq.~\eqref{eq2:thm:Tensor Mean Inequality Under Linear Transform}.
$\hfill \Box$

\bibliographystyle{IEEETran}
\bibliography{AHType_2Argu_Bib}

\end{document}